\documentclass[10pt,reqno]{amsart}
\usepackage[T1]{fontenc}
\usepackage{lmodern,microtype}
\usepackage{amsmath,amssymb,amsthm,mathtools}
\usepackage{booktabs,longtable,ltcaption,array,enumitem}

\usepackage{xcolor,hyperref}
\usepackage[capitalise,noabbrev]{cleveref}
\hypersetup{colorlinks=true,linkcolor=blue!50!black,citecolor=green!40!black,urlcolor=blue!60!black,
	pdftitle={Recurrence formulas for Taylor expansions at CM points on genus-one and genus-two Shimura curves},
	pdfauthor={Jingxuan Geng}}
\allowdisplaybreaks
\numberwithin{equation}{section}
\newtheorem{theorem}{Theorem}[section]
\newtheorem{proposition}[theorem]{Proposition}
\newtheorem{lemma}[theorem]{Lemma}
\newtheorem{corollary}[theorem]{Corollary}
\theoremstyle{definition}
\theoremstyle{remark}

\newcommand{\Q}{\mathbb{Q}}
\newcommand{\C}{\mathbb{C}}
\newcommand{\cR}{\mathcal R}
\newcommand{\cO}{\mathcal O}

\newcommand{\Sch}[2]{\left\{#1,#2\right\}}

\title[Recurrences on genus-one and genus-two Shimura curves]
{Recurrence formulas for Taylor expansions at CM points on genus-one and genus-two Shimura curves}
\author{Jingxuan Geng}
\begin{document}
	\begin{abstract}
		We establish a general differential equation relating the expansion of a weight-two eigenform $f$ on a Shimura curve to the Schwarzian derivative $\mathcal{R}$ associated with $f$. When the quotient $X/W$ has genus zero, we introduce a Hauptmodul $u$ and obtain a coupled system of differential equations involving $f$, $\mathcal{R}$, and $u$. Geometric information, including local constraints on the poles of $\mathcal{R}$, together with the initial coefficients of the expansion of $f$, determine all the parameters in the system. We apply this method to all twelve weight-two eigenforms on Shimura curves $X_0^N(1)$ of genus one or two, obtaining explicit recurrence formulas for their Taylor expansions at selected CM points, conditional on the numerically recognized initial coefficients.  As an application, we use the coefficients obtained from the recurrence to compute the Shimura uniformization $\Phi: \operatorname{Jac}(X)\to E$ to high precision and identify a large family of CM points on the associated elliptic curves $E$.

	\end{abstract}
	\maketitle
	
	\section{Introduction}
	\label{sec:intro}
	
	\subsection{Explicit computation on Shimura curves}
	
	Let \(X_N=X_0^N(1)\) be the Shimura curve of trivial level associated with the indefinite quaternion algebra over \(\Q\) of discriminant \(N\). Explicit equations, uniformization maps, and evaluations of modular forms are fundamental tools for studying the arithmetic of these Shimura curves.
	The absence of cusps on $X_N$ makes these computations particularly challenging, since the usual methods for classical modular curves rely on Fourier expansions at cusps.
	
	When $X_N$ has small genus, one can often obtain explicit models by combining coverings with arithmetic information about CM points. For example, Gonz\'alez--Rotger \cite{GR2004,GR} determine models of genus-one and genus-two Shimura curves by studying their elliptic Jacobians and elliptic quotients. Guo--Yang \cite{GuoYang2017} give equations for hyperelliptic Shimura curves by constructing modular functions with prescribed CM divisors as Borcherds forms. These results yield \cref{tab:models}, which lists the models
	\(y^2=F_N(x)\) for Shimura curves $X_N$ of genus one or two, together with the eigendifferentials used in this paper.
	
	\begin{table}[!htbp]
		\centering\small
		\renewcommand{\arraystretch}{1.55}
		\setlength{\tabcolsep}{5pt}
		\caption{Models and eigendifferentials for all \(X_N\)
			of genus one from Gonz\'alez--Rotger \cite[Theorem~3.4]{GR} and two from Guo--Yang \cite[Appendix~A, Table~2]{GuoYang2017}.}
		\label{tab:models}
		\begin{tabular}{@{}lll@{}}
			\toprule
			\(N\)&Model&Eigendifferential \(\omega\)\\
			\midrule
			\(14\)&\(y^2=-x^4+13x^2-128\)&\(dx/y\)\\
			\(15\)&\(y^2=-3x^4-82x^2-27\)&\(dx/y\)\\
			\(21\)&\(y^2=-7x^4+94x^2-343\)&\(dx/y\)\\
			\(33\)&\(y^2=-3x^4-10x^2-243\)&\(dx/y\)\\
			\(34\)&\(y^2=-3x^4+26x^3-53x^2-26x-3\)&\(dx/y\)\\
			\(46\)&\(y^2=-x^4+45x^2-512\)&\(dx/y\)\\
			\midrule
			\(26\)&\(y^2=-2x^6+19x^4-24x^2-169\)&\(\begin{aligned}\omega_a&=x\,dx/y,\\\omega_b&=dx/y.\end{aligned}\)\\
			\(38\)&\(y^2=-16x^6-59x^4-82x^2-19\)&\(\begin{aligned}\omega_a&=x\,dx/y,\\\omega_b&=dx/y.\end{aligned}\)\\
			\(58\)&\(y^2=-2x^6-78x^4-862x^2-1682\)&\(\begin{aligned}\omega_a&=x\,dx/y,\\\omega_b&=dx/y.\end{aligned}\)\\
			\bottomrule
		\end{tabular}
	\end{table}

	The Schwarz\-ian derivative and the associated differential equations are useful tools for explicit computations on Shimura curves. Elkies \cite{Elkies} introduced these tools to study the local inverse $z$ of the covering map $p: \mathcal{H} \to X_N$. Using ramification data and Atkin--Lehner symmetries, Elkies computed covering maps and their local inverses numerically and algebraically and constructed equations for the curves. Tu \cite{Tu} determines Schwarzian equations for a family of genus-zero Shimura quotients and obtains the Schwarzian derivative for many curves of small genus. Yang \cite{Yang2013} uses solutions of Schwarzian differential equations to compute a basis for the space of modular forms $S_k(\Gamma)$ and to compute Hecke operators on genus-zero Shimura curves $X_N=\mathcal{H}/\Gamma$.
	
	The analytic computation of modular forms on a compact Shimura curve is another fundamental problem. In the absence of cusp expansions, one can instead study the Taylor expansion at an interior point $\tau \in X_N$, such as a CM point, as proposed by Zagier \cite[Section~6.3]{Zagier}. See also Voight--Willis \cite[Section~2]{VW}. For an eigenform \(f\) on $X_N$ and a CM point \(\tau\in\mathcal H\), define the \(\tau\)-centered coordinate
	\begin{equation}\label{eq:w}
		w=\frac{z-\tau}{z-\bar\tau},
	\end{equation}
	and consider the expansion of $f$ at $w=0$ given by
	\begin{equation}\label{eq:f(w)}
		f(z(w))=(1-w)^k\sum_{n\ge0}b_nw^n,
	\end{equation}
	where the map \(z\mapsto w\) sends \(\mathcal H\) to the unit disk \(\mathcal D\) and
	\(\tau\) to \(0\). 
	
	When \(b_0b_1\ne0\), it is convenient to use the normalized series $A(q)$ obtained by factoring out the period factor $\Theta=\frac{b_1}{b_0}$ and setting $q=\Theta w$, as in \cite[Section~2]{VW}.
	\begin{equation}\label{eq:A(q)}
		A(q)=\frac{f(\frac{q}{\Theta})}{(1-w)^kb_0}=\sum_{n\ge0}c_n\frac{q^n}{n!},
		\qquad c_n=n!\frac{b_n}{b_0\Theta^n}.
	\end{equation}
	
	The reason for working with $A(q)$ is that CM theory gives \(c_n\in\overline{\Q}\). For all the forms computed in this paper, we have \(c_n\in \Q\). If we normalize $f$ by \(f(\tau)=b_0=1\), we recover the original modular form $f(w)$ from $A(q)$ by
	\begin{equation}\label{eq:f(thetaw)}
		f(w)=(1-w)^kA(\Theta w)
		=(1-w)^k\sum_{n\ge0}\frac{c_n}{n!}(\Theta w)^n.
	\end{equation}
	We review this normalization in Section \ref{sec:disk-coefficients}.
	
	Several methods are available for computing the expansions $f(w)$ and $A(q)$ numerically. Voight--Willis
	\cite[Section~3]{VW} compute the coefficients using modularity relations and Cauchy's integral formula. Nelson \cite[Sections~1.5 and~3.4]{Nelson} uses the explicit Shimizu theta
	correspondence to evaluate $b_n$ as certain lattice sums. An approach based on the formulas of
	Bertolini--Darmon--Prasanna can compute $b_n^2$ by relating it to central \(L\)-values of Rankin--Selberg products of $f$ and CM theta series \cite[Appendix~A.2, Project~B]{Prasanna}.
	
	These numerical methods work well for a small number of initial coefficients $b_n$, but computing many more coefficients is costly. For example, for the weight-two eigenform $14a$, the method above readily gives $b_n$ up to $n=100$, whereas computing $b_n$ up to $n=10000$ is out of reach. This limits the precision attainable from the computed expansion of $f$.
	
	Voight--Willis \cite[Remark~2.6]{VW} raise the problem of finding explicit recursions for CM expansions by studying the differential structure of modular forms on Shimura curves. Such recursions can easily determine all subsequent coefficients exactly from finitely many initial coefficients. Zagier \cite[Proposition~28, p.~88]{Zagier} combines Ramanujan's differential identities with modular differential operators to obtain a polynomial recursion for the expansion of the Eisenstein series $E_4$ at the CM point of discriminant \(-4\). Baba--Granath \cite{BG} adapt this approach to modular forms on the Shimura curve $X_{6}$: they use the Schwarzian differential equation of a modular function to establish a differential system for generators of the modular-form ring, from which they derive explicit recurrences for the expansions.
	
	\subsection{Main result}
	Given an eigenform $f$ on a Shimura curve $X_N$, we first present a differential equation satisfied by the normalized form $A(q)$ in \eqref{eq:A(q)} and the Schwarzian derivative $\cR$ associated to $f$, which is a meromorphic function on $X_N$ whose definition will be recalled in Section \ref{sec:Schwarzian}. 
	\begin{theorem}\label{thm:recurrence}
		Let \(f\) be a weight-two eigenform on $X_N$ and $\tau \in \mathcal{H}$.
		Suppose the expansion \eqref{eq:f(w)} of $f$ at $w_\tau$ has \(b_0b_1\ne0\), and let \(A(q)\) be its normalized series. Then $A(q)$ and the associated Schwarzian derivative $\cR(q)$ satisfy the differential equation
		\begin{equation}\label{eq:main-recurrence}
			AA''-\frac32(A')^2+L\cR(q)A^4=0.
		\end{equation}
		Here $L$ is a constant, and the derivatives are taken with respect to $q$.
	\end{theorem}
	The proof of Theorem~\ref{thm:recurrence} is a direct application of the chain rule for the Schwarzian derivative and will be presented in Section \ref{sec:Schwarzian}.  The immediate consequence of Theorem~\ref{thm:recurrence} is that the Schwarzian derivative $\cR$ essentially determines the expansion of the eigenform: one can completely determine the expansion of $A(q)$ and thus $f(w)$ if one knows the expansion of the Schwarzian derivative $\cR(q)$ in terms of $q$ and the constant $L$. This is never easy, as the expansion of Schwarzian derivative $\cR(q)$ is in general hard to compute as well. 
	
	However, when the genus $g(X_N)$ is small, we may use the existing models to describe the Schwarzian derivative $\cR$. To be more precise, as the Schwarzian derivative $\cR$ is invariant under the action of the Atkin--Lehner group $W_N$ (see Proposition \ref{prop:descent}), when the quotient curve $X_N/W$ has genus zero, we can use the Hauptmodul $u$ of $X_N/W$ to express $\cR=\cR(u)$ and form a closed coupled recurrence system involving $A(q)$, $u(q)$ and $\cR(u(q))$. 
	\begin{theorem}\label{thm:coupled-recurrence}
		Suppose $X_N/W$ has genus zero and let $u=u(q)$ be a Hauptmodul of $X_N/W$. By \cref{prop:descent}, $\cR=\cR(u) \in \mathbb{C}(u)$. There exists a function $H(u)\in \mathbb{C}(u)$ such that 
		\begin{equation}\label{eq:coupled-recurrence}
			\boxed{\begin{aligned}
					AA''-\frac32(A')^2+L\cR(u)A^4&=0,\\
					u''-\frac{A'}A u'-\frac L2 A^2H'(u)&=0.
			\end{aligned}}
		\end{equation}
		Here, $A$, $\cR$, and $L$ are as in Theorem \ref{thm:recurrence}. The derivatives of $A$ and $u$ are taken with respect to $q$, while $H'(u)$ denotes the derivative with respect to $u$.
	\end{theorem}
	
	From \eqref{eq:A(q)} we have $c_0=c_1=1$. Thus the Schwarzian derivative $\cR(u)$, the constant $L$, the function $H(u)$, and the initial values $u(0)$ and $u'(0)$ uniquely determine the coupled recurrence. More precisely, if we write \(u(q)=\sum_{m\ge0}u_m\frac{q^m}{m!},\) then \(c_0=c_1=1\) and \(u_0,u_1\) determine every \(c_n\) and \(u_n\) recursively by
	\[
	\begin{aligned}
		c_{n+2}
		&=
		n!\,[q^n]\left(
		\frac{3(A')^2}{2A}-L\cR(u)A^3
		\right),\\
		u_{n+2}
		&=
		n!\,[q^n]\left(
		\frac{A'}{A}u'
		+\frac L2 A^2H'(u)
		\right).
	\end{aligned}
	\]
	Here $[q^n]F(q)$ denotes the coefficient of $q^n$ in the expansion of $F(q)$.

	Conversely, sufficiently many initial coefficients, computed numerically by the methods of \cite{Nelson,Prasanna,VW}, can help us recover the unknown parameters, conditional on the algebraic recognition of the $c_n$. We combine geometric data from the literature, including the Schwarzian derivatives determined for several levels in \cite{Tu}, with initial coefficients $c_n$ computed numerically using \cite{Nelson,Prasanna,VW}. This determines all the parameters in the coupled recurrence system of Theorem~\ref{thm:coupled-recurrence} for all weight-two eigenforms on genus-one and genus-two Shimura curves, as stated in \cref{thm:main}.

	\begin{theorem}\label{thm:main}
		For each eigenform $f$ in \cref{tab:main-data}, let $A(q)$ be its normalized expansion at a CM point of discriminant $D$ on $X_N$. For each \(f\), assume the values of the initial coefficients $c_n$ listed in \cref{tab:input-data} and \cref{tab:cm-expansions}. Then \cref{tab:main-data} gives the function $u$, the Schwarzian derivative $\cR(u)$, the constant $L$, the function $H(u)$, and the initial values $u(0)$ and $u'(0)$. These data uniquely determine the differential equation \eqref{eq:coupled-recurrence} and all subsequent coefficients $c_n$ of $A(q)$ and $u_n$ of $u(q)$.
		
	\end{theorem}
	\begingroup
	\small
	\setlength{\tabcolsep}{0.6pt}
	\renewcommand{\arraystretch}{1.35}
	\setlength{\extrarowheight}{3pt}
	\setlength{\LTcapwidth}{\linewidth}
	\setlength{\LTleft}{\fill}
	\setlength{\LTright}{\fill}
	\begin{longtable}{@{}
			>{\centering\arraybackslash}m{0.8cm}
			>{\centering\arraybackslash}m{0.75cm}
			>{\centering\arraybackslash}m{1.3cm}
			>{\centering\arraybackslash}m{0.8cm}
			>{\centering\arraybackslash}m{1.7cm}
			>{\centering\arraybackslash}m{1.7cm}
			>{\centering\arraybackslash}m{2.7cm}
			>{\centering\arraybackslash}m{2.6cm}
			@{}}
		\caption{CM centers and recurrence data for all twelve eigenforms.}
		\label{tab:main-data}\\
		\toprule
		Form&\(D\)&\(u\)&\(u_0\)&\(u_1=u'(0)\)&\(L\)&\(\cR(u)\)&\(H(u)\)\\
		\midrule
		\endfirsthead
		\multicolumn{8}{c}{\tablename~\thetable\ (continued)}\\
		\toprule
		Form&\(D\)&\(u\)&\(u_0\)&\(u_1=u'(0)\)&\(L\)&\(\cR(u)\)&\(H(u)\)\\
		\midrule
		\endhead
		\midrule
		\multicolumn{8}{r}{\footnotesize Continued on the next page}\\
		\endfoot
		\bottomrule
		\endlastfoot
		\(\displaystyle 14a\) & \(\displaystyle -11\) & \(\displaystyle x^2\) & \(\displaystyle 16\) & \(\displaystyle -\frac{352}{3}\) & \(\displaystyle -\frac{11}{9}\) & \(\displaystyle \frac38(4-u)\) & \(\displaystyle \begin{gathered}4u(-u^2+13u\\[-2pt]{}-128)\end{gathered}\) \\
		\specialrule{0.3pt}{3pt}{3pt}
		\(\displaystyle 15a\) & \(\displaystyle -7\) & \(\displaystyle x^2\) & \(\displaystyle 9\) & \(\displaystyle -63\) & \(\displaystyle -\frac7{64}\) & \(\displaystyle -4-\frac43u\) & \(\displaystyle \begin{gathered}4u(-3u^2-82u\\[-2pt]{}-27)\end{gathered}\) \\
		\specialrule{0.3pt}{3pt}{3pt}
		\(\displaystyle 21a\) & \(\displaystyle -43\) & \(\displaystyle x^2\) & \(\displaystyle \frac{441}{25}\) & \(\displaystyle -\frac{9632}{375}\) & \(\displaystyle -\frac{43}{3969}\) & \(\displaystyle -\frac{\begin{gathered}16(17u^2+392u\\[-2pt]{}-2401)\end{gathered}}{(u-49)^2}\) & \(\displaystyle \begin{gathered}4u(-7u^2+94u\\[-2pt]{}-343)\end{gathered}\) \\
		\specialrule{0.3pt}{3pt}{3pt}
		\(\displaystyle 33a\) & \(\displaystyle -67\) & \(\displaystyle x^2\) & \(\displaystyle \frac9{25}\) & \(\displaystyle -\frac{19296}{2875}\) & \(\displaystyle -\frac{67}{529}\) & \(\displaystyle -\frac{\begin{gathered}4(u^3+49u^2\\[-2pt]{}+19u+27)\end{gathered}}{u(u-1)^2}\) & \(\displaystyle \begin{gathered}4u(-3u^2-10u\\[-2pt]{}-243)\end{gathered}\) \\
		\specialrule{0.3pt}{3pt}{3pt}
		\(\displaystyle 34a\) & \(\displaystyle -11\) & \(\displaystyle x-x^{-1}\) & \(\displaystyle -\frac{15}{4}\) & \(\displaystyle -\frac{3179}{80}\) & \(\displaystyle -\frac{11}{25}\) & \(\displaystyle \frac{\begin{gathered}279u^2+3104u\\[-2pt]{}-39264\end{gathered}}{2(15u-16)^2}\) & \(\displaystyle \begin{gathered}(u^2+4)\\[-2pt]{}\cdot(-3u^2+26u\\[-2pt]{}-59)\end{gathered}\) \\
		\specialrule{0.3pt}{3pt}{3pt}
		\(\displaystyle 46a\) & \(\displaystyle -163\) & \(\displaystyle x^2\) & \(\displaystyle \frac{3136}{81}\) & \(\displaystyle -\frac{15188992}{636255}\) & \(\displaystyle -\frac{844992}{61701025}\) & \(\displaystyle -\frac{\begin{gathered}3u^3+52u^2\\[-2pt]{}+13824u\\[-2pt]{}-311296\end{gathered}}{8(u-64)^2}\) & \(\displaystyle \begin{gathered}4u(-u^2+45u\\[-2pt]{}-512)\end{gathered}\) \\
		\specialrule{0.5pt}{3pt}{3pt}
		\(\displaystyle 26a\) & \(\displaystyle -11\) & \(\displaystyle x^2\) & \(\displaystyle 1\) & \(\displaystyle \frac{44}{21}\) & \(\displaystyle -\frac{11}{1764}\) & \(\displaystyle \frac92-\frac9u-\frac{507}{2u^2}\) & \(\displaystyle \begin{aligned}&4\bigl(-2u^3+19u^2\\&\quad-24u-169\bigr)\end{aligned}\) \\
		\specialrule{0.3pt}{3pt}{3pt}
		\(\displaystyle 26b\) & \(\displaystyle -11\) & \(\displaystyle x^{-2}\) & \(\displaystyle 1\) & \(\displaystyle 44\) & \(\displaystyle -\frac{11}{4}\) & \(\displaystyle 3+\frac{14}u-\frac3{u^2}\) & \(\displaystyle \begin{aligned}&4\bigl(-169u^3\\&\quad-24u^2\\&\quad+19u-2\bigr)\end{aligned}\) \\
		\specialrule{0.3pt}{3pt}{3pt}
		\(\displaystyle 38a\) & \(\displaystyle -11\) & \(\displaystyle x^2\) & \(\displaystyle 1\) & \(\displaystyle \frac{44}{9}\) & \(\displaystyle -\frac{11}{324}\) & \(\displaystyle \begin{gathered}-6u-\frac{45}2\\[2pt]{}-\frac{57}u-\frac{57}{2u^2}\end{gathered}\) & \(\displaystyle \begin{aligned}&4\bigl(-16u^3\\&\quad-59u^2\\&\quad-82u-19\bigr)\end{aligned}\) \\
		\specialrule{0.3pt}{3pt}{3pt}
		\(\displaystyle 38b\) & \(\displaystyle -11\) & \(\displaystyle x^{-2}\) & \(\displaystyle 1\) & \(\displaystyle \frac{44}{13}\) & \(\displaystyle -\frac{11}{676}\) & \(\displaystyle -16-\frac{52}u-\frac{30}{u^2}\) & \(\displaystyle \begin{aligned}&4\bigl(-19u^3\\&\quad-82u^2\\&\quad-59u-16\bigr)\end{aligned}\) \\
		\specialrule{0.3pt}{3pt}{3pt}
		\(\displaystyle 58a\) & \(\displaystyle -11\) & \(\displaystyle x^2\) & \(\displaystyle -11\) & \(\displaystyle -\frac{176}{5}\) & \(\displaystyle \frac{121}{400}\) & \(\displaystyle \begin{gathered}5-\frac{144358}{243u}\\[2pt]{}-\frac{2523}{u^2}\\[2pt]{}-\frac{229376}{243(u+27)}\\[2pt]{}+\frac{65536}{9(u+27)^2}\end{gathered}\) & \(\displaystyle \begin{aligned}&4\bigl(-2u^3-78u^2\\&\quad-862u-1682\bigr)\end{aligned}\) \\
		\specialrule{0.3pt}{3pt}{3pt}
		\(\displaystyle 58b\) & \(\displaystyle -11\) & \(\displaystyle x^{-2}\) & \(\displaystyle -\frac1{11}\) & \(\displaystyle -\frac{16}{33}\) & \(\displaystyle -\frac{11}{144}\) & \(\displaystyle \begin{gathered}-\frac{39625}{243}\\[2pt]{}-\frac{34}u-\frac3{u^2}\\[2pt]{}-\frac{163840}{243(27u+1)}\\[2pt]{}-\frac{65536}{243(27u+1)^2}\end{gathered}\) & \(\displaystyle \begin{aligned}&4\bigl(-1682u^3\\&\quad-862u^2\\&\quad-78u-2\bigr)\end{aligned}\) \\
	\end{longtable}
	\endgroup

	Consequently, we determine the Taylor expansions of all such eigenforms $f$, given in Appendix \ref{sec:cm-expansion}, conditional on the algebraic recognition of $\Theta_N$ and the initial $c_n$.
	
	\begin{corollary}[CM expansions]\label{cor:cm-expansions}
		Conditional on the initial values of $c_n$ and the values of $\Theta_N$ listed in \cref{tab:cm-expansions}, the forms normalized by \(f_N(\tau_N)=1\)
		have the CM expansions in \cref{tab:cm-expansions}.
		All subsequent coefficients $c_n$ are determined by the recurrence \eqref{eq:coupled-recurrence}.
	\end{corollary}
	
	Another immediate corollary is that the coefficients $c_n$, known a priori to lie in $\bar{\mathbb{Q}}$, are rational. Moreover, $c_n\in\mathbb{Z}[\frac{1}{S}]$ for some $S\in \mathbb{Z}$.
	\begin{proposition}\label{prop:Z[1/S]}
		Conditional on the initial values of $c_n$ in Table~\ref{tab:cm-expansions}, for each form $f$ in Table~\ref{tab:main-data} we have
		\[
		c_n\in\mathbb{Z}[1/S_f]\qquad\text{for all } n\ge 0,
		\]
		where the minimal integer $S_f$ is given by
		\[
		\setlength{\arraycolsep}{3.2pt}
		\begin{array}{c|cccccccccccc}
			f   & 14a & 15a & 21a & 33a & 34a & 46a & 26a & 26b & 38a & 38b & 58a & 58b\\ \hline
			S_f & 3 & 2 & 21 & 23 & 5 & 15710 & 21 & 1 & 3 & 13 & 5 & 3
		\end{array}
		\]
	\end{proposition}
	The proof, which will be given in Section \ref{sec:Z[1/S]}, uses a change of coordinate \(v=\phi(u)\) for which the coupled system \cref{eq:coupled-recurrence} is \(\mathbb{Z}[1/S_f]\)-integral in $v$.
	
	We use a combination of Sage~\cite{SageMath} and Magma~\cite{Magma} to carry out all computations.
	
	\subsection{Outline}
	Section \ref{sec:background} reviews the definition of modular forms on Shimura curves and the normalization of their expansions at CM points. Section \ref{sec:Schwarzian} defines the Schwarzian derivative, studies it using Atkin--Lehner descent and Riemann--Roch, and proves Theorem~\ref{thm:recurrence} and Theorem~\ref{thm:coupled-recurrence}. Section \ref{sec:determine-parameters} determines the remaining parameters in \cref{tab:cm-expansions} from existing geometric data and the assumed coefficients $c_n$ and completes the proofs of \cref{thm:main} and \cref{prop:Z[1/S]}. Appendix \ref{sec:cm-expansion} lists the assumed initial coefficients $c_n$ and gives the expansions of the eigenforms.
	Appendix~\ref{app:u-values} certifies CM values of the Hauptmodul $u$ used in the paper.
	
	\subsection{Acknowledgement}
	The author would like to thank John Voight for suggesting the research problem. The author appreciates Michael Lipnowski and Kartik Prasanna for valuable discussions and comments on this work.
	\section{Modular forms and CM expansions on Shimura curves}
	\label{sec:background}
	
	\subsection{Shimura curves, modular forms, and differentials}
	\label{sec:shimura-curves}
	\label{sec:differentials}
	
	Let \(B_N\) be the indefinite quaternion algebra over \(\Q\) of
	discriminant \(N>1\), and let \(\cO\subset B_N\) be a maximal order.
	Fix a real splitting with algebraic matrix entries, and write
	\(\cO^1\) for the group of elements of reduced norm one. The associated
	Shimura curve \(X_N\) has complex uniformization
	\cite[Sections~2.2--2.3]{Nelson}
	\[
	X_N(\C)=\mathcal H/\Gamma,\qquad
	\Gamma=\cO^1/\{\pm1\}\subset\mathrm{PSL}_2(\mathbb{R}).
	\]
	For an even integer \(k\ge2\), a modular form of weight \(k\) on \(X_N\)
	is a holomorphic function \(f:\mathcal H\to\C\) satisfying
	\[
	f\!\left(\frac{az+b}{cz+d}\right)=(cz+d)^k f(z)
	\qquad
	\left(\begin{pmatrix}a&b\\c&d\end{pmatrix}\in\cO^1\right).
	\]
	We denote this space by \(S_k(\Gamma)\).
	
	A holomorphic differential on \(X_N\) is a section of
	\(\Omega_{X_N}^1\). Under the uniformization
	\(p:\mathcal H\to\mathcal H/\Gamma=X_N\), weight-two
	modular forms realize these differentials analytically through the following Hecke-equivariant isomorphism
	\cite[Section~2.3]{Nelson}:
	\[
	S_2(\Gamma)\simeq H^0(X_N,\Omega_{X_N}^1),\qquad
	f\longmapsto\eta_f,\qquad p^*\eta_f=f(z)\,dz.
	\]
	
	Fix a nonzero weight-two eigenform \(f\). Choose a nonzero holomorphic
	differential \(\omega\) in the line spanned by \(\eta_f\), and write
	\(\eta_f=c_f\omega\), with \(c_f\in\C^\times\).
	One can define a local primitive coordinate \(t\) of \(\omega\) by
	\begin{equation}\label{eq:local_t}
		dt=\omega=\frac{1}{c_f}f(z)dz
	\end{equation}
	It is a local coordinate wherever \(\omega\ne0\). For the genus-one and genus-two curves, we take \(t\) from the eigendifferentials listed in \cref{tab:models}.
	For a genus-one curve, \(\omega=dx/y\) is a nowhere-vanishing differential, and \(t\) is a local Abel coordinate \cite[Theorem~21.10]{Forster}. For a genus-two curve, $\omega_a$ has simple zeros at the two points above $x=0$, while $\omega_b$ has simple zeros at the two points above $x=\infty$.

	\subsection{Normalization of CM expansions}
	\label{sec:disk-coefficients}
	We recall the standard Taylor expansion of a modular form at an interior point of a Shimura curve \cite[Section~2]{VW}. As in \eqref{eq:f(w)}, we choose $\tau \in \mathcal{H}$ and use \[w=\frac{z-\tau}{z-\bar\tau}\] as our parameter. For \(f\in S_k(\Gamma)\), the disk expansion is
	\[
	f(z(w))=(1-w)^kB(w),\qquad B(w)=\sum_{n\ge0}b_nw^n.
	\]
	The inverse coordinate change gives
	\begin{equation}\label{eq:disk}
		z=\frac{\tau-\bar\tau w}{1-w},\qquad
		f(z)\,dz=(\tau-\bar\tau)(1-w)^{k-2}B(w)\,dw.
	\end{equation}
	
	The coefficients \(b_n\) are values of Shimura--Maass derivatives.
	With
	\[
	\delta_k=\frac{1}{2\pi i}
	\left(\frac{\partial}{\partial z}+\frac{k}{z-\bar z}\right),
	\qquad
	\delta_k^{(n)}=\delta_{k+2n-2}\cdots\delta_{k+2}\delta_k,\qquad
	\delta_k^{(0)}=\mathrm{id},
	\]
	the Taylor identity \cite[Lemma~2.2]{VW} gives
	\[
	b_n=\frac{(-4\pi\operatorname{Im}\tau)^n}{n!}
	\delta_k^{(n)}f(\tau).
	\]
	Moreover, it is natural to choose a CM point \(\tau\) as our expansion center, as the $b_n$ for CM points are algebraic up to a normalization factor factoring out the following Chowla--Selberg period defined by
	\cite[Section~6.3, equation~(97)]{Zagier}
	\[
	\Omega_D=\frac{1}{\sqrt{2\pi|D_K|}}
	\left(\prod_{j=1}^{|D_K|-1}
	\Gamma\!\left(\frac{j}{|D_K|}\right)^{\chi_{D_K}(j)}\right)^{w_K/(4h_K)},
	\]
	where $D_K$ is the fundamental discriminant of the CM point $\tau$, \(h_K\) is the class number of \(K\), \(w_K=|\cO_K^\times|\), and \(\chi_{D_K}(j)=\left(\frac{D_K}{j}\right)\) is the Kronecker symbol. The coefficients $b_n$ satisfy the algebraicity relation \cite[Theorem~2.5]{VW}
	\[
	\frac{b_n}{\pi^n\Omega_D^{2n}}
	=\frac{(-4\operatorname{Im}\tau)^n}{n!}
	\frac{\delta_k^{(n)}f(\tau)}{\Omega_D^{2n}}
	\in\overline{\Q}.
	\]
	When \(b_0b_1\ne0\), we may construct a normalized form $A(q)$ with algebraic coefficients as in
	\cite[Section~2, equations~(5)--(7)]{VW} by extracting the period factor $\Theta$ given below.
	\[
	\begin{aligned}
		\Theta
		&= \frac{b_1}{b_0} \in \pi\Omega_D^2\overline{\Q}^{\times}, \\[4pt]
		q
		&= \Theta w, \\[4pt]
		A(q)
		&= \frac{f\!\left(\frac{q}{\Theta}\right)}{(1-w)^k b_0}
		= \sum_{n\ge 0} c_n \frac{q^n}{n!},\\
		c_n &= n!\frac{b_n}{b_0\Theta^n}.
	\end{aligned}
	\]
	Thus \(c_0=c_1=1\) and \(c_n\in\overline{\Q}\). In our
	normalization \(b_0=1\), this expresses the original modular form $f$ as
	\[
	f(z)=(1-w)^kA(\Theta w)
	=(1-w)^k\sum_{n\ge0}\frac{c_n}{n!}(\Theta w)^n.
	\]
	
	Specializing to \(k=2\), at a center with \(b_0b_1\ne0\), the differential \(\omega\) does not
	vanish and \(t\) is a local coordinate. Combining
	\(f(z)\,dz=c_f\,dt\) with \eqref{eq:disk} and
	\(q=\Theta w\) gives the comparison with the algebraic differential:
	\begin{equation}\label{eq:dt=lamdaA(q)}
		dt=\lambda A(q)\,dq,\qquad
		\lambda=\frac{(\tau-\bar\tau)b_0}{c_f\Theta}.
	\end{equation}
	Thus integrating $A(q)$ gives the local coordinate $t$ up to the factor $\lambda$.

	\section{The Schwarzian derivative}
	\label{sec:Schwarzian}
	
	We define the Schwarzian derivative, study its properties by Akin-Lehner descent and Riemann-Roch, and prove \cref{thm:recurrence} and \cref{thm:coupled-recurrence}.
	
	\subsection{The Schwarzian derivative}
	
	For a function \(g(x)\), its Schwarzian derivative is
	\cite[Section~2, p.~455]{Tu}
	\[
	\Sch{g}{x}=\frac{g'''}{g'}-\frac32\left(\frac{g''}{g'}\right)^2.
	\]
	It measures the failure of \(g\) to be locally M\"obius and vanishes
	precisely when \(g\) is a local M\"obius transformation.
	It satisfies the M\"obius invariance
	\begin{equation}\label{eq:schwarzian-mobius}
		\Sch{M\circ g}{x}=\Sch{g}{x},
	\end{equation}
	where \(M\) is a M\"obius transformation. The Schwarzian also satisfies the chain rule
	\cite[Section~2.2, equations~(13)--(14)]{Elkies}
	\begin{equation}\label{eq:schwarzian-chain}
		\Sch{g\circ h}{x}
		=\left.\Sch{g}{s}\right|_{s=h(x)}
		\left(\frac{dh}{dx}\right)^2+\Sch{h}{x}.
	\end{equation}
	Here the derivatives in \(\Sch{g}{s}\) are taken with respect to \(s=h(x)\).
	
	Elkies \cite{Elkies} studied the Schwarzian derivative 	
	\begin{equation}\label{eq:R-definition}
		\cR=\{z,t\},
	\end{equation} 
	of the local inverse $z\colon X_N \to \mathcal{H}$ of the uniformization \(p: \mathcal{H} \to X_N \), with respect to a local coordinate such as $t$ in \eqref{eq:local_t}. The branches of \(z\) differ by elements of \(\Gamma\), so M\"obius invariance \eqref{eq:schwarzian-mobius} ensures that the Schwarzian derivative is well defined away
	from the elliptic points of $X_N$ and the zeros of \(\omega\). Moreover, \(\cR\) extends meromorphically to $X_N$, with prescribed pole orders and leading polar coefficients described by the following lemma. 
	\begin{lemma}
		\label{lem:local-schwarzian}
		Set \(Z=\operatorname{div}(\omega)\) and \(D_{\mathrm{ell}}=\sum_{P\text{ elliptic}}P\),
		using the same symbols for their supports. The Schwarzian \(\cR=\{z,t\}\) of
		\eqref{eq:R-definition} extends meromorphically to \(X_N\) with poles at \(Z\cup D_{\mathrm{ell}}\).
		
		At each pole \(P\), let \(m=m_P=\operatorname{ord}_P(\omega)\ge0\) and let \(e=e_P\)
		be its elliptic order.
		In a local coordinate \(v\) centered at \(P\), write \(\omega=v^m(a+bv+O(v^2))\,dv\)
		with \(a\ne0\). Then the pole of \(\cR\) at \(P\) has order \(2m+2\), with
		\begin{equation}\label{eq:local-pole}
			\cR(v)=\frac{(m+1)^2-e^{-2}}{2a^2}\,v^{-2m-2}+O(v^{-2m-1}).
		\end{equation}
		For \(P\in Z\setminus D_{\mathrm{ell}}\), \(\omega\) determines coefficients at \(v^{-2m-1}\) as well, with
		\[
		\cR(v)=\frac{m(m+2)}{2a^2}\,v^{-2m-2}
		-\frac{m(m+1)b}{a^3}\,v^{-2m-1}+O(v^{-2m}).
		\]
		For \(P\in D_{\mathrm{ell}}\setminus Z\), the coordinate \(v=t-t(P)\) gives
		\begin{equation}\label{eq:local-pole-elliptic}
			\cR(v)=\kappa_e v^{-2}+O(v^{-1}),\qquad \kappa_e=\frac{1-e^{-2}}{2}.
		\end{equation}
	\end{lemma}
	
	\begin{proof}
		Write \(dt=\omega=h(v)\,dv\). On a punctured neighborhood of \(P\), the Schwarzian chain rule gives
		\begin{equation}\label{eq:chain-R-h}
			\cR(v)
			=
			\frac{\{z,v\}-\{t,v\}}{h(v)^2}.
		\end{equation}
		
		At a nonelliptic point, the local inverse uniformization \(z=z(v)\)
		is holomorphic with nonzero derivative, so \(\{z,v\}\) is holomorphic.
		At an elliptic point of order \(e\), local uniformization gives the expansion (see \cite[\S2.2]{Elkies})
		\[
		\{z,v\}
		=
		\frac{1-e^{-2}}{2v^2}+O(v^{-1});
		\]

		From the expansion of \(h=av^m+bv^{m+1}+O(v^{m+2})\), we have
		\[
		\frac{h'}{h}
		=
		\frac{m}{v}+\frac{b}{a}+O(v).
		\]
		Consequently,
		\[
		\{t,v\}
		=
		\left(\frac{h'}{h}\right)'
		-\frac12\left(\frac{h'}{h}\right)^2
		=
		-\frac{m(m+2)}{2v^2}
		-\frac{mb}{av}
		+O(1).
		\]
		Together with \eqref{eq:chain-R-h} and \(h(v)^{-2}=a^{-2}v^{-2m}\left(1-\frac{2b}{a}v+O(v^2)\right),\)
		this yields \eqref{eq:local-pole}. When \(e=1\), the holomorphy
		of \(\{z,v\}\) gives the stated two-term expansion.
		
		These formulas give the meromorphic extension across \(P\).
		The leading coefficient in \eqref{eq:local-pole} is nonzero
		whenever \(m>0\) or \(e>1\), proving the asserted pole order.
		When \(m=0\) and \(e=1\), the same identity shows that \(\cR\)
		is holomorphic. Finally, if \(\omega(P)\ne0\), then \(t-t(P)\)
		is a regular local coordinate, and choosing this coordinate gives
		\eqref{eq:local-pole-elliptic}.
	\end{proof}

	We now prove the main recurrence equation \eqref{eq:main-recurrence} in \cref{thm:recurrence}. The simple but key observation is that since \(q=\Theta(z-\tau)/(z-\bar\tau)\), \(z\) is M\"obius in \(q\). Thus \(\Sch{z}{q}=0\). Then a direct application of the chain rule of Schwarzian derivatives connecting \(\Sch{z}{q}\) and \(\Sch{z}{t}\) gives a proof of \eqref{eq:main-recurrence}.
	\begin{proof}[Proof of \cref{thm:recurrence}]
		By \eqref{eq:dt=lamdaA(q)}, we have $\frac{dt}{dq}=\lambda A(q)$.
		By the chain rule \eqref{eq:schwarzian-chain}, we have
		\begin{equation}\label{eq:schwarzian-recurrence-chain}
			\begin{aligned}
				0=\{z,q\}
				&=\{z,t\}\left(\frac{dt}{dq}\right)^2+\{t,q\}\\
				&=\cR(t(q))\lambda^2A(q)^2+\{t,q\}\\
				&=\cR(t(q))\lambda^2A(q)^2+\frac{A''}{A}
				-\frac32\left(\frac{A'}A\right)^2.
			\end{aligned}
		\end{equation}
		Multiplying the equation by $A^2$ and setting $L=\lambda^2$ gives \eqref{eq:main-recurrence}.
	\end{proof}
	
	One can see that we can determine the recurrence if we have a good understanding of \(\cR\), which we will study next using Atkin--Lehner descent and Riemann--Roch.
	
	\subsection{Atkin--Lehner descent}
	The Atkin--Lehner group $W_N$ acts by M\"obius transformations on $z$. The following proposition shows that $\cR$ is invariant under $W_N$.
	\begin{proposition}[Atkin--Lehner descent]\label{prop:descent}
		Let \(\omega\ne0\) be a differential associated with an eigenform on \(X_N\).
		For \(\cR=\{z,t\}\), with \(dt=\omega\), we have
		\[
		\cR\in\C(X_N)^{W_N}.
		\]
		If \(X_N/W_N\) has genus zero with Hauptmodul \(u\), then
		\[
		\cR,\ H=\left(\frac{du}{dt}\right)^2
		\in\C(X_N)^{W_N}=\C(u).
		\]
		In particular, this holds for all eigenforms with the coordinates
		\(u\) listed in \cref{tab:main-data}.
	\end{proposition}
	
	\begin{proof}
		Let \(w\in W_N\). The eigenform differential satisfies
		\[
		w^*\omega=\varepsilon_w\omega,\qquad \varepsilon_w\in\{\pm1\}.
		\]
		Fix a point \(P\) away from the elliptic points and the zeros of
		\(\omega\). The same holds for \(wP\), since \(w^*\omega=\varepsilon_w\omega\).
		Choose local Abel coordinates
		\(t\) near \(P\) and \(s\) near \(wP\), each obtained by integrating
		\(\omega\). The coordinate expression of \(w\) is
		\[
		h=s\circ w\circ t^{-1},\qquad
		h(t)=\varepsilon_w t+c_w.
		\]
		In particular, \(dh/dt=\varepsilon_w\) and \(\{h,t\}=0\).
		Let \(z_P(t)\) and \(z_{wP}(s)\) be local inverse uniformizations
		expressed in these coordinates. Choose the branches compatibly with
		a M\"obius lift \(M_w\) of \(w\), so that
		\(z_{wP}\circ h=M_w\circ z_P\).
		M\"obius invariance \eqref{eq:schwarzian-mobius} gives
		\[
		\{z_{wP}\circ h,t\}
		=\{M_w\circ z_P,t\}=\{z_P,t\}.
		\]
		On the other hand, the chain rule \eqref{eq:schwarzian-chain} gives
		\[
		\begin{aligned}
			\{z_{wP}\circ h,t\}
			&=\left.\{z_{wP},s\}\right|_{s=h(t)}
			\left(\frac{dh}{dt}\right)^2+\{h,t\}\\
			&=\left.\{z_{wP},s\}\right|_{s=h(t)}.
		\end{aligned}
		\]
		Evaluation at \(t=t(P)\) gives
		\(\cR(P)=\cR(wP)\);
		the equality extends across the excluded points by meromorphicity.
		
		If \(u\) is a Hauptmodul on \(X_N/W_N\), differentiating \(u\circ w=u\)
		gives \(\varepsilon_w(du/ds)(wP)=(du/dt)(P)\). Squaring gives
		\(H=\left(\frac{du}{dt}\right)^2\in\C(X_N)^{W_N}\) as $\varepsilon_w^2=1$. 
		
		For all curves in \cref{tab:models}, the action of \(W_N\) is given in \cite[Theorem~3.4]{GR} for genus one and
		\cite[Proposition~2.1 and  Theorem~3.1]{GR2004} for genus two. It is generated by
		\[
		\iota(x,y)=(x,-y),\qquad
		\sigma(x,y)=
		\begin{cases}
			(-x,y),& N\ne34,\\
			(-1/x,-y/x^2),&N=34.
		\end{cases}
		\]
		These generators fix
		\(u=x^2\), or \(u=x-x^{-1}\) for \(N=34\); since \(u\) has degree
		four and \(W_N\) has order four, its invariant field is exactly \(\C(u)\).
	\end{proof}
	
	We now use \cref{prop:descent} to deduce the coupled differential equations \eqref{eq:coupled-recurrence}.
	
	\begin{proof}[Proof of \cref{thm:coupled-recurrence}]
		Write $u_t=\frac{du}{dt}$ and $u_{tt}=\frac{d^2u}{dt^2}$. The identity \(u_t^2=H(u)\) gives \(u_{tt}=H'(u)/2\).   Using \(t_q=\lambda A\) and
		\(L=\lambda^2\), we obtain 
		\[
		u'=\frac{du}{dt}\frac{dt}{dq}=u_t \lambda A,
		\]  
		and 
		\begin{equation}\label{eq:coordinate-recurrence}
			\begin{aligned}
				u''&= u_{tt}(\lambda A)^2+u_t\lambda A'\\
				&=\frac{H'(u)}{2}LA^2+\frac{A'}A u'.
			\end{aligned}
		\end{equation}
		Thus we have proved \eqref{eq:coupled-recurrence}.
	\end{proof}

	The function \(H(u)=u_t^2\) can be computed
	directly for all forms in \cref{tab:families} from the explicit models and differentials in \cref{tab:models}. For genus-one curves, except when $N=34$, the model has the form \(y^2=P_N(x^2)\) with quadratic \(P_N\) and \(u=x^2\). We have \(dt=dx/y\), whence
	\[
	\frac{du}{dt}=2xy,
	\]
	and
	\begin{equation}\label{eq:H-computation}
		H_N(u)=\left(\frac{du}{dt}\right)^2
		=4x^2y^2=4u(a_Nu^2+b_Nu+c_N).
	\end{equation}
	
	For \(N=34\), the coordinate \(u=x-x^{-1}\) and the quartic model of $X_{34}$ gives
	\[
	\begin{aligned}
		H_{34}(u)=\left(\frac{du}{dt}\right)^2&=(x+x^{-1})^2\frac{y^2}{x^2}\\
		&=(x+x^{-1})^2\frac{-3x^4+26x^3-53x^2-26x-3}{x^2}\\
		&=(u^2+4)(-3u^2+26u-59).
	\end{aligned}
	\]
	
	For the two eigenforms \(\omega_a\) and \(\omega_b\) on a genus-two curve with model \(y^2=P_N(x^2)\) with cubic \(P_N\), the same calculation gives
	\[
	\frac{du_a}{dt_a}=2y,\qquad
	\frac{du_b}{dt_b}=-\frac{2y}{x^3},
	\]
	and therefore
	\begin{equation}\label{eq:H-genus-two}
		H_{N,a}(u)=4P_N(u),\qquad H_{N,b}(u)=4u^3P_N(1/u).
	\end{equation}
	This proves the formulas in the $H_N(u)$ column of \cref{tab:main-data}.

	\subsection{Riemann--Roch and the form of
		\texorpdfstring{\(\cR(u)\)}{R(u)}}
	\label{sec:riemann-roch}
	
	We now determine the finite-dimensional space in which the Schwarzian derivative \(\cR(u)\) must lie by Riemann--Roch after imposing the local conditions stated in \cref{lem:local-schwarzian}.
	In particular, \(\cR\) is holomorphic outside \(Z\cup D_{\mathrm{ell}}\). At a point
	\(P\in D_{\mathrm{ell}}\setminus Z\), the coordinate
	\(v=t-t(P)\) gives
	\[
	\cR(v)=\kappa_{e_P}v^{-2}+O(v^{-1}),
	\qquad
	\kappa_{e_P}=\frac{1-e_P^{-2}}2,
	\]
	as in \eqref{eq:local-pole-elliptic}.
	For the genus-one forms, \(Z=0\), so this describes all
	poles of \(\cR\).
	
	For the genus-two forms in \cref{tab:models}, \(Z\)
	consists of two simple zeros: those of
	\(\omega_a=x\,dx/y\) lie above \(x=0\), while those of
	\(\omega_b=dx/y\) lie at infinity. At these points, take
	\(v=x\) and \(v=1/x\), respectively, and put
	\(\eta=y/x^3\) near infinity. The even sextic models of the genus-two curves give
	\[
	\omega_a
	=
	\left(\frac{v}{y(P)}+O(v^3)\right)dv,
	\qquad
	\omega_b
	=
	\left(-\frac{v}{\eta(P)}+O(v^3)\right)dv.
	\]
	At these zeros, the coefficients in
	\cref{lem:local-schwarzian} satisfy \(m=1\) and \(b=0\). If we write \(P_N(s)=a_Ns^3+b_Ns^2+c_Ns+d_N\), then
	\(a^{-2}=d_N\) for \(\omega_a\) and \(a^{-2}=a_N\) for \(\omega_b\).
	So all parameters are explicit and the lemma gives
	\[
	\cR(v)=
	\begin{cases}
		\displaystyle
		\frac{3}{2a^2}\,v^{-4}+O(v^{-2}),
		& P\in Z\setminus D_{\mathrm{ell}},\\[6pt]
		\displaystyle
		\frac{4-e_P^{-2}}{2a^2}\,v^{-4}+O(v^{-3}),
		& P\in Z\cap D_{\mathrm{ell}}.
	\end{cases}
	\]
	In particular, on \(Z\setminus D_{\mathrm{ell}}\),
	\(\cR\) has the prescribed coefficient of \(v^{-4}\) and
	has zero coefficient of \(v^{-3}\).
	
	Now suppose a global meromorphic function \(\cR_0\) on \(X_N\) satisfies all the local conditions above. Thus
	\(\cR_0\) has no poles outside \(Z\cup D_{\mathrm{ell}}\)
	and has the same pole orders and prescribed polar coefficients as the Schwarzian derivative. It follows that
	\[
	\cR-\cR_0=
	\begin{cases}
		O(v^{-1}),
		& P\in D_{\mathrm{ell}}\setminus Z,\\
		O(v^{-2}),
		& P\in Z\setminus D_{\mathrm{ell}},\\
		O(v^{-3}),
		& P\in Z\cap D_{\mathrm{ell}},\\
		O(1),
		& P\notin Z\cup D_{\mathrm{ell}}.
	\end{cases}
	\]
	For a divisor \(E\), write
	\(L_{X_N}(E)=H^0(X_N,\cO_{X_N}(E))\).
	Since \(Z\) is reduced in our cases, we therefore obtain
	\[
	\cR\in\cR_0+L_{X_N}(2Z+D_{\mathrm{ell}}).
	\]

	Put \(s=\deg D_{\mathrm{ell}}\). For the genus-one forms,
	\(Z=0\) and \(s>0\), while for the genus-two forms,
	\(\deg Z=2\). In both cases,
	\(\deg(2Z+D_{\mathrm{ell}})>2g(X_N)-2\), so
	Riemann--Roch gives
	\[
	\dim L_{X_N}(2Z+D_{\mathrm{ell}})
	=
	\begin{cases}
		s,   & g(X_N)=1,\\
		s+3, & g(X_N)=2.
	\end{cases}
	\]
	
	The prescribed local conditions are stable under \(W_N\).
	By averaging a reference function over \(W_N\), we may therefore
	choose \(\cR_0\) to be \(W_N\)-invariant.
	Proposition~\ref{prop:descent} then gives
	\[
	\cR\in
	\cR_0+L_{X_N}(2Z+D_{\mathrm{ell}})^{W_N}.
	\]
	It is very helpful to describe the latter invariant space by pullback from the quotient
	\[
	\pi:X_N\longrightarrow X_N/W_N=\mathbf P^1_u.
	\]
	Define
	\[
	\Delta_\omega
	=
	\sum_{a\in\mathcal U}
	\left\lfloor
	\frac{2m_a+\mathbf1_{\{e_a>1\}}}{r_a}
	\right\rfloor[a].
	\]
	Here \(\mathcal U=\pi(Z\cup D_{\mathrm{ell}})\), and
	for \(a\in\mathcal U\), we write \(m_a\), \(e_a\), and \(r_a\)
	for the vanishing order of \(\omega\), the elliptic order on
	\(X_N\), and the ramification index of \(\pi\), respectively,
	at any point above \(a\). Since
	\(\operatorname{ord}_P(\pi^*F)
	=r_a\operatorname{ord}_a(F)\) for \(P\in\pi^{-1}(a)\),
	we have
	\[
	L_{X_N}(2Z+D_{\mathrm{ell}})^{W_N}
	=
	\pi^*L_{\mathbf P^1_u}(\Delta_\omega).
	\]
	Since \(W_N\)
	is an abelian \(2\)-group, \(r_a=1\) when no nontrivial
	Atkin--Lehner involution fixes a point above \(a\), and
	\(r_a=2\) when exactly one does.

	We can compute the fixed polar coefficients downstairs,
	using
	\[
	H(u)=u_t^2=\left(\frac{du}{\omega}\right)^2.
	\]
	There are three types of finite pole fibers in our examples.
	
	First suppose that \(u=a\) is an unramified elliptic fiber
	with \(m_a=0\). Then \(H(a)\ne0\), and
	\[
	u-a=u_t(P)(t-t(P))+O((t-t(P))^2).
	\]
	Using \eqref{eq:local-pole-elliptic}, we obtain
	\begin{equation}\label{eq:unramified-pole}
		\cR(u)
		=
		\frac{\kappa_{e_a}H(a)}{(u-a)^2}
		+\frac{B_a}{u-a}
		+O(1).
	\end{equation}
	The double-pole coefficient is fixed, while the coefficient
	\(B_a\) remains to be determined.
	
	Next suppose that \(u=b\) is a ramified elliptic fiber
	with \(m_b=0\). We may write
	\[
	u-b=c(t-t(P))^2+O((t-t(P))^3),
	\qquad c\ne0.
	\]
	Comparing the leading terms in \(u_t^2=H(u)\) gives
	\[
	H(b)=0,
	\qquad
	4c^2=H'(b)c,
	\qquad
	c=\frac{H'(b)}4.
	\]
	The double pole upstairs descends to a simple pole, whose
	coefficient is \(c\kappa_{e_b}\). Thus
	\begin{equation}\label{eq:ramified-pole}
		\cR(u)
		=
		\frac{\kappa_{e_b}H'(b)}{4(u-b)}
		+O(1).
	\end{equation}
	In particular, this fiber contributes nothing to \(\Delta_\omega\), since
	\(\lfloor1/2\rfloor=0\).
	
	Finally, suppose that \(u=a\) is a zero fiber.
	In our genus-two examples, \(m_a=1\) and \(r_a=2\).
	Choose a regular local coordinate \(v\) with \(u-a=v^2\),
	and write
	\[
	\omega=(\alpha v+O(v^2))\,dv,
	\qquad \alpha\ne0.
	\]
	Then \(H(a)=4/\alpha^2\ne0\).
	By \cref{lem:local-schwarzian}, the fourth-order pole
	upstairs descends to a double pole with coefficient
	\[
	\frac{4-e_a^{-2}}{2\alpha^2}
	=
	\kappa_{2e_a}H(a),
	\qquad
	\kappa_{2e_a}=\frac{4-e_a^{-2}}8.
	\]
	Consequently,
	\begin{equation}\label{eq:zero-pole}
		\cR(u)
		=
		\frac{\kappa_{2e_a}H(a)}{(u-a)^2}
		+\frac{B_a}{u-a}
		+O(1).
	\end{equation}
	Whether or not the zero is elliptic, its contribution to
	\(\Delta_\omega\) is
	\[
	\left\lfloor
	\frac{2+\mathbf1_{\{e_a>1\}}}{2}
	\right\rfloor[a]=[a].
	\]
	
	Now consider the pole fibers at infinity. This occurs for the forms
	\(14a\), \(15a\), \(46a\), and \(38a\).
	In each case \(u=x^2\), \(m_\infty=0\), and \(r_\infty=2\).
	Put \(v=1/x\). For the genus-one quartics,
	\[
	\frac{dv}{dt}=-\frac{y}{x^2},
	\]
	while for the genus-two form \(\omega_a=x\,dx/y\),
	\[
	\frac{dv}{dt}=-\frac{y}{x^3}.
	\]
	Hence, in both cases,
	\[
	\left(\frac{dv}{dt}(P)\right)^2=a_N,
	\]
	where \(a_N\) is the leading coefficient of the corresponding $P_N$. Since \(u=v^{-2}\), the leading
	coefficient in \eqref{eq:local-pole-elliptic} gives
	\begin{equation}\label{eq:infinite-pole}
		\cR(u)=\kappa_{e_\infty}a_Nu+O(1).
	\end{equation}

	Let \(\mathcal U_{\mathrm{fin}}
	=\mathcal U\setminus\{\infty\}\).
	Combining the preceding calculations gives
	\begin{equation}\label{eq:general-accessory-family}
		\begin{aligned}
			\cR(u)={}&C+\kappa_{e_\infty}a_Nu
			+\sum_{\substack{a\in\mathcal U_{\mathrm{fin}}\\
					(m_a,r_a)=(0,1)}}
			\left(
			\frac{\kappa_{e_a}H(a)}{(u-a)^2}
			+\frac{B_a}{u-a}
			\right)\\
			&+\sum_{\substack{a\in\mathcal U_{\mathrm{fin}}\\
					(m_a,r_a)=(0,2)}}
			\frac{\kappa_{e_a}H'(a)}{4(u-a)}
			+\sum_{\substack{a\in\mathcal U_{\mathrm{fin}}\\
					(m_a,r_a)=(1,2)}}
			\left(
			\frac{\kappa_{2e_a}H(a)}{(u-a)^2}
			+\frac{B_a}{u-a}
			\right).
		\end{aligned}
	\end{equation}

	For example, for the form \(46a\), the elliptic fibers are
	\(u=64\), of order three and ramification index one,
	and \(u=\infty\), of order two and ramification index two. In the notation \((a,e_a,r_a,m_a)\) in \cref{tab:families} below, the corresponding tuples are \((\infty,2,2,0)\) and \((64,3,1,0)\). All the $u$-values used here are certified in Appendix \ref{app:u-values}.
	Here
	\[
	H_{46}(u)=4u(-u^2+45u-512),
	\qquad
	H_{46}(64)=-442368.
	\]
	Equation \eqref{eq:general-accessory-family} then gives
	\begin{equation}\label{eq:general-accessory-family-46}
		\cR_{46}(u)
		=
		C+\frac{B}{u-64}
		-\frac{196608}{(u-64)^2}
		-\frac{3u}{8}.
	\end{equation}
	Indeed, \(\kappa_3H_{46}(64)=-196608\), and the
	coefficient at infinity is \(\kappa_2(-1)=-3/8\).
	Moreover,
	\[
	\Delta_\omega=[64],
	\qquad
	L_{X_{46}}(D_{\mathrm{ell}})^{W_{46}}
	=
	\pi^*L_{\mathbf P^1_u}([64]).
	\]
	Similar computations in the accompanying computational certificate in Sage, give all the families
	in \cref{tab:families}.

	\begingroup
	\small
	\setlength{\tabcolsep}{3pt}
	\renewcommand{\arraystretch}{1.5}
	\setlength{\LTcapwidth}{\linewidth}
	\setlength{\LTleft}{\fill}
	\setlength{\LTright}{\fill}
	\begin{longtable}{@{}
			>{\centering\arraybackslash}p{0.8cm}
			>{\centering\arraybackslash}p{3cm}
			>{\raggedright\arraybackslash}p{7.8cm}
			>{\centering\arraybackslash}p{0.4cm}@{}}
		\caption{Families containing $\cR$. A tuple
			\((a,e_a,r_a,m_a)\) records the fiber over \(u=a\), its elliptic order \(e_a\) on \(X_N\), ramification index \(r_a\in\{1,2\}\) of \(\pi:X_N\to X_N/W_N\), and vanishing order \(m_a\) of \(\omega\). Here \(d=1+\operatorname{deg}(\Delta_\omega)\) is the number of remaining parameters.}
		\label{tab:families}\\
		\toprule
		Form&Pole fibers \((a,e_a,r_a,m_a)\)&Family for \(\cR(u)\)&\(d\)\\
		\midrule
		\endfirsthead
		\multicolumn{4}{c}{\tablename~\thetable\ (continued)}\\
		\toprule
		Form&Pole fibers \((a,e_a,r_a,m_a)\)&Family for \(\cR(u)\)&\(d\)\\
		\midrule
		\endhead
		\midrule
		\multicolumn{4}{r}{\footnotesize Continued on the next page}\\
		\endfoot
		\bottomrule
		\endlastfoot
		\(14a\) & \((\infty,2,2,0)\) & \(\displaystyle C-\frac{3u}{8}\) & \(1\)\\
		\specialrule{0.3pt}{5pt}{5pt}
		\(15a\) & \((\infty,3,2,0)\) & \(\displaystyle C-\frac{4u}{3}\) & \(1\)\\
		\specialrule{0.3pt}{5pt}{5pt}
		\(21a\) & \((49,2,1,0)\) & \(\displaystyle C+\frac{B}{u-49}-\frac{921984}{(u-49)^2}\) & \(2\)\\
		\specialrule{0.3pt}{5pt}{5pt}
		\(33a\) & \(\begin{gathered}(0,3,2,0),\\(1,2,1,0)\end{gathered}\) & \(\displaystyle C+\frac{B}{u-1}-\frac{384}{(u-1)^2}-\frac{108}{u}\) & \(2\)\\
		\specialrule{0.3pt}{5pt}{5pt}
		\(34a\) & \((16/15,3,1,0)\)  & \(\displaystyle C+\frac{B}{u-16/15}-\frac{1336336}{16875(u-16/15)^2} \) & \(2\)\\
		\specialrule{0.3pt}{5pt}{5pt}
		\(46a\) & \(\begin{gathered}(\infty,2,2,0),\\(64,3,1,0)\end{gathered}\) & \(\displaystyle C+\frac{B}{u-64}-\frac{196608}{(u-64)^2}-\frac{3u}{8}\) & \(2\)\\
		\specialrule{0.5pt}{5pt}{5pt}
		\(26a\) & \((0,1,2,1)\) & \(\displaystyle C+\frac{B}{u}-\frac{507}{2u^2}\) & \(2\)\\
		\specialrule{0.3pt}{5pt}{5pt}
		\(26b\) & \((0,1,2,1)\) & \(\displaystyle C+\frac{B}{u}-\frac{3}{u^2}\) & \(2\)\\
		\specialrule{0.3pt}{5pt}{5pt}
		\(38a\) & \(\begin{gathered}(0,1,2,1),\\(\infty,2,2,0)\end{gathered}\) & \(\displaystyle C+\frac{B}{u}-\frac{57}{2u^2}-6u\) & \(2\)\\
		\specialrule{0.3pt}{5pt}{5pt}
		\(38b\) & \((0,2,2,1)\) & \(\displaystyle C+\frac{B}{u}-\frac{30}{u^2}\) & \(2\)\\
		\specialrule{0.3pt}{5pt}{5pt}
		\(58a\) & \(\begin{gathered}(0,1,2,1),\\(-27,3,1,0)\end{gathered}\) & \(\displaystyle C+\frac{B_0}{u}+\frac{B_1}{u+27}-\frac{2523}{u^2}+\frac{65536}{9(u+27)^2}\) & \(3\)\\
		\specialrule{0.3pt}{5pt}{5pt}
		\(58b\) & \(\begin{gathered}(0,1,2,1),\\(-1/27,3,1,0)\end{gathered}\) & \(\displaystyle C+\frac{B_0}{u}+\frac{B_1}{u+1/27}-\frac{3}{u^2}-\frac{65536}{177147(u+1/27)^2}\) & \(3\)\\
		
	\end{longtable}
	\endgroup
	
	\section{Determination of \texorpdfstring{$\cR$, $L$, and $u_1$}{R, L, and u\_1}}
	\label{sec:determine-parameters}
	Having determined $u$, $H(u)$, and the possible forms of $\cR(u)$, we determine the remaining data $\cR(u)$, $L$, and $u_1$ in \cref{tab:main-data}, conditional on the values of the $c_n$. These data uniquely determine the subsequent coefficients $c_n$. Conversely, sufficiently many coefficients $c_n$ beyond $c_1=c_0=1$ determine the parameters uniquely.
	
	The determination is divided into two cases: one group has existing literature on the Schwarzian $\cR(u)$, for which fewer coefficients $c_n$ are needed. The other group has a free parameter in its $\cR(u)$, for which more coefficients $c_n$ are needed.

	\subsection{Forms with known \texorpdfstring{\(\cR(u)\)}{R(u)}}
	\label{sec:14}
	We first consider \(14a\), \(15a\), and \(21a\), for which Tu \cite[Section~5.1]{Tu} obtains closely related Schwarzian equations. A change of variable and the Schwarzian chain rule give the derivative \(\cR(u)\) in our coordinates.
	
	Once \(\cR(u)\) is known, we can determine the remaining parameters $L$ and $u_1$ by rewriting \eqref{eq:main-recurrence} as
	\begin{equation}\label{eq:reconstruction}
		L\cR(u(q))=\frac{\frac32(A')^2-AA''}{A^4}.
	\end{equation}
	Since \(c_0=c_1=1\), differentiation gives
	\begin{equation}\label{eq:q-jets}
		\begin{aligned}
			L\cR(u_0)&=\frac32-c_2,\\
			\left.\frac{d}{dq}\bigl(L\cR(u(q))\bigr)\right|_{q=0}
			&=6c_2-c_3-6,
		\end{aligned}
	\end{equation}
	These two identities then determine
	\begin{equation}\label{eq:scale-branch-recovery}
		L=\frac{3/2-c_2}{\cR(u_0)},\qquad
		u_1=\frac{6c_2-c_3-6}{L\cR'(u_0)},
	\end{equation}
	Thus $c_2$ determines $L$, and the additional value of $c_3$ determines $u_1$.

	For \(14a\), let \(dt=dx/y\). We have
	\begin{equation}\label{eq:14-data}
		y^2=-x^4+13x^2-128,\qquad u=x^2,\qquad
		H_{14}(u)=u_t^2=4u(-u^2+13u-128).
	\end{equation}
	
	In the model $v^2=-16s^4+13s^2-8$ of $X_{14}$ given by Elkies \cite{Elkies}, Tu \cite[Section~5.1]{Tu} chooses $T=-s^2$ and computes
	\[
	Q_{14}(T)=\frac12\{z,T\}
	=
	\frac{192+440T+43T^2+1036T^3+960T^4}
	{16T^2(8+13T+16T^2)^2}.
	\]
	A change of variable \(x=4s,\ y=4v\) will recover our model and identify
	\[
	T=\phi(u)=-\frac{u}{16}
	\]

	Table \ref{tab:families} gives \(\cR_{14}(u)=C-\frac38u\).
	The change \(x=4s,\ y=4v\) preserves \(dx/y=ds/v\). Then We can recover $\cR(u)=\Sch{z}{t}$ from Tu's $\Sch{z}{T}$ using $T=\phi(u)$ and the Schwarzian chain rule \eqref{eq:schwarzian-chain}.

	\begin{equation}\label{eq:quotient-conversion}
		\begin{aligned}
			\cR(u)&=\{z,t\}\\
			&=\{z,u\}u_t^2+\{u,t\}\\
			&=\bigl(\{z,T\}\phi'(u)^2+\{\phi,u\}\bigr)H(u)
			+\frac{u_{ttt}}{u_t}
			-\frac32\left(\frac{u_{tt}}{u_t}\right)^2\\
			&=2H(u)Q(\phi(u))\phi'(u)^2+H(u)\{\phi,u\}
			+\frac12H''(u)-\frac{3H'(u)^2}{8H(u)}.
		\end{aligned}
	\end{equation}
	
	In the last step we used \(u_{tt}=H'(u)/2\) and \(u_{ttt}=H''(u)u_t/2\), obtained by differentiating \(u_t^2=H(u)\). These identities extend meromorphically across the zeros of \(u_t\).
	
	Since \(\phi(u)=-u/16\), we have \(\{\phi,u\}=0\) and
	\(\phi'(u)^2=1/256\). Substituting the displayed \(H_{14}\) and \(Q_{14}\) and simplifying 
	gives
	\begin{equation}\label{eq:14-R}
		\begin{aligned}
			\cR_{14}(u)
			&=\frac{H_{14}(u)}{128}Q_{14}\!\left(-\frac{u}{16}\right)
			+\frac12H_{14}''(u)
			-\frac{3H_{14}'(u)^2}{8H_{14}(u)}\\
			&=\frac38(4-u).
		\end{aligned}
	\end{equation}
	This agrees with our computation in \cref{tab:families}.

	With \(\cR_{14}\) known, we determine $L$ and $u_1$ from \eqref{eq:scale-branch-recovery} using the numerical values of $c_2$ and $c_3$ given by
	\begin{equation}\label{eq:14-seeds}
		c_2=-4,\qquad c_3=214/9.
	\end{equation}
	We obtain 
	\begin{equation}\label{eq:14-initial}
		L=\frac{3/2-c_2}{\cR_{14}(u_0)}=-\frac{11}{9},\quad
		\quad u_1=\frac{6c_2-c_3-6}{L\cR_{14}'(u_0)}=-\frac{352}{3}.
	\end{equation}
	This completes the determination of the parameters for the form $14a$. 
	
	The parameters for $15a$ and $21a$ can be determined in exactly the same way. For $15a$, 
	Tu takes \(T=-3x^2=-3u\) and obtains
	\[
	Q_{15}(T)=
	\frac{177147-244944T+244242T^2-3680T^3+35T^4}
	{144T^2(1-T)^2(81-T)^2}.
	\]
	Using \(\phi(u)=-3u\) in \eqref{eq:quotient-conversion} gives \(\cR_{15}=-4-\frac43u\). As in \eqref{eq:14-initial}, the coefficients $c_2=-\frac14$ and $c_3=\frac{27}{16}$ determine  
	\begin{equation}\label{eq:15-initial}
		L=-\frac{7}{64},\quad
		\quad u_1=-63.
	\end{equation}
	
	For $21a$, Tu's coordinate $T=u$ coincides with ours, and the formula is
	\begin{equation}\label{eq:21-Q}
		\begin{aligned}
			Q_{21}(u)&=\frac{21\left(\begin{gathered}40353607-17647350u+3561369u^2-477652u^3\\{}+31833u^4-630u^5+7u^6\end{gathered}\right)}{16u^2(49-u)^2(343-94u+7u^2)^2},\\
		\end{aligned}
	\end{equation}
	
	Applying \eqref{eq:quotient-conversion} with $\phi(u)=u$ gives
	\begin{equation}\label{eq:21-R}
		\cR_{21}(u)=-272-\frac{32928}{u-49}-\frac{921984}{(u-49)^2}.
	\end{equation}
	Then the values \(c_2=-905/3969\) and \(c_3=-1318/27783\) determine \begin{equation}\label{eq:21-scale}
		L=-\frac{43}{3969},\qquad
		u_1=-\frac{9632}{375}.
	\end{equation}
	
	\subsection{Forms without known \texorpdfstring{\(\cR(u)\)}{R(u)}}
	\label{sec:high}
	We now recover the parameters from additional coefficients $c_n$ when $\cR(u)$ is only partially known. Taking higher derivatives as in \eqref{eq:q-jets} gives
	\begin{equation}\label{eq:q-jets-3}
		\begin{aligned}
			\left.\frac{d^2}{dq^2}\bigl(L\cR(u(q))\bigr)\right|_{q=0}
			&=6c_2^2+9c_3-c_4-42c_2+30,\\
			\left.\frac{d^3}{dq^3}\bigl(L\cR(u(q))\bigr)\right|_{q=0}
			&=21c_2c_3-108c_2^2+330c_2-78c_3+12c_4-c_5-180.
		\end{aligned}
	\end{equation}.
	
	By \cref{tab:families}, for each remaining form except $58a$ and $58b$, the Schwarzian $\cR$ belongs to the following two-parameter family.
	\[
	\cR(u)=C+\frac{B}{u-a}+\Phi(u),
	\]
	where \(\Phi\) denotes the known part.
	The coupled equation \eqref{eq:coupled-recurrence} expresses \(u_2,u_3\) in terms of \(u_1\) and \(L\):
	\begin{equation}\label{eq:fit-u2-u3}
		\begin{aligned}
			u_2&=u_1+\frac L2H'(u_0),\\
			u_3&=c_2u_1+\frac{3L}{2}H'(u_0)+\frac L2H''(u_0)u_1.
		\end{aligned}
	\end{equation}
	The chain rule expresses the derivatives in \eqref{eq:q-jets} and \eqref{eq:q-jets-3} as
	\begin{equation}\label{eq:q-matching}
		\begin{aligned}
			\left.\frac{d}{dq}\bigl(L\cR(u(q))\bigr)\right|_{q=0}
			&=L\cR'(u_0)u_1,\\
			\left.\frac{d^2}{dq^2}\bigl(L\cR(u(q))\bigr)\right|_{q=0}
			&=L\bigl(\cR''(u_0)u_1^2+\cR'(u_0)u_2\bigr),\\
			\left.\frac{d^3}{dq^3}\bigl(L\cR(u(q))\bigr)\right|_{q=0}
			&=L\bigl(\cR'''(u_0)u_1^3
			+3\cR''(u_0)u_1u_2+\cR'(u_0)u_3\bigr).
		\end{aligned}
	\end{equation}
	Here primes on \(H\) and \(\cR\) denote derivatives with respect to \(u\).
	The identity \(u_1^2=LH(u_0)\), \eqref{eq:q-jets}, and \eqref{eq:q-matching} express \(L,B,C\) directly in terms of \(u_1\) and the known
	coefficients \(c_2,c_3\):
	\begin{equation}\label{eq:fit-elimination}
		\begin{aligned}
			L&=\frac{u_1^2}{H(u_0)},\\
			B&=(u_0-a)^2\left(\Phi'(u_0)
			-\frac{H(u_0)(6c_2-c_3-6)}{u_1^3}\right),\\
			C&=\frac{H(u_0)(3/2-c_2)}{u_1^2}
			-\frac{B}{u_0-a}-\Phi(u_0).
		\end{aligned}
	\end{equation}
	Thus it suffices to determine \(u_1\).
	
	Observe that the left-hand sides of \eqref{eq:q-matching} are known values if we know values up to $c_5$, while on the right-hand sides every term can now be expressed in terms of $u_1$. Substituting these expressions would yield polynomial equations in $u_1$. Solving them would determine $u_1$.
	
	Let us take the form \(26b\) as an example.
	By \cref{tab:families},
	\[
	\cR_{26,b}(u)=C+\frac{B}{u}-\frac{3}{u^2}.
	\]
	Assume the values of $c_2 \cdots c_5$ in \cref{tab:cm-expansions}.
	Then, by \eqref{eq:fit-u2-u3} and the known $H_{26}$,
	\[
	\begin{aligned}
		u_2&=u_1+\frac{67}{44}u_1^2,\\
		u_3&=40u_1+\frac{201}{44}u_1^2+\frac{531}{176}u_1^3.
	\end{aligned}
	\]
	
	For \(26b\), \eqref{eq:fit-elimination} gives
	\begin{equation}\label{eq:26b-elimination}
		\begin{aligned}
			L&=-\frac{u_1^2}{704},\\
			B&=6+\frac{681472}{u_1^3},\\
			C&=-3+\frac{27104}{u_1^2}-\frac{681472}{u_1^3}.
		\end{aligned}
	\end{equation}
	
	Then we can solve for $u_1$ by the following equations
	\[
	0=E_j(u_1):=\left.\frac{d^j}{dq^j}
	\left(L\cR_{26,b}(u(q))-
	\frac{\frac32(A')^2-AA''}{A^4}\right)\right|_{q=0},
	\qquad j=2,3,
	\]
	Note that every term can be expressed in terms of $u_1$,
	and we compute
	\begin{equation}\label{eq:26b-fit-polynomials}
		\begin{aligned}
			E_2(u_1)&=-11616-462u_1
			+\frac{3}{352}u_1^4,\\[3pt]
			E_3(u_1)&=2201232-1386u_1-\frac{231}{2}u_1^2
			+\frac{9}{352}u_1^4-\frac{189}{15488}u_1^5.
		\end{aligned}
	\end{equation}
	The unique common solution over \(\C\) of $E_2(u_1)=E_3(u_1)=0$ is
	\[
	u_1=44,
	\]
	and consequently \(L=-11/4\),
	\(B=14\), and \(C=3\).
	
	The same argument applies to the remaining cases. For every remaining form except $58a$, the values of $c_2 \dots c_5$ determine all the parameters through $E_2$ and $E_3$. For $58a$, we also consider $u_4$ and $E_4$, which requires the additional value $c_6$. The Sage certificate records the polynomials $E_i$ and the computation for each form. This completes the proof of \cref{thm:main}.
	
	\subsection{Integrality of \texorpdfstring{$c_n$}{c\_n}}\label{sec:Z[1/S]}
	To prove \cref{prop:Z[1/S]}, it is easier to work with the following equation
	\[
	A'''=6\frac{A'A''}{A}-6\frac{(A')^3}{A^2}
	-L\cR'(u)u'A^3.
	\]
	This is obtained by differentiation of the first equation (\ref{eq:coupled-recurrence}). The advantage is all coeffcients here are integral and there is no $\frac{3}{2}$ factor causing spurious denominators.
	
	Then we choose an appropriate fractional linear coordinate \(v=\phi(u)\) with \(\phi(u_0)=0\), so that
	\(v'(0)=\phi'(u_0)u_1\). Using \((u')^2=LA^2H(u)\), the system (\ref{eq:coupled-recurrence}) in $v$ is 
	\begin{equation}\label{eq:integral-system}
		\begin{aligned}
			A'''&=6\frac{A'A''}{A}-6\frac{(A')^3}{A^2}-F(v)v'A^3,\\
			v''&=\frac{A'}A v'+G(v)A^2,
		\end{aligned}
	\end{equation}
	where
	\begin{equation}\label{eq:integral-FG}
		\begin{aligned}
			F(v)&=\left.\frac{L\cR'(u)}{\phi'(u)}\right|_{u=\phi^{-1}(v)},\\
			G(v)&=\left.\frac L2\bigl(2\phi''(u)H(u)+\phi'(u)H'(u)\bigr)
			\right|_{u=\phi^{-1}(v)}.
		\end{aligned}
	\end{equation}
	A suitable choice of $v$ with ideal $P$ and $Q$ will make everything in (\ref{eq:integral-system}) $S$-integral, as the following lemma shows.
	\begin{lemma}[Integrality criterion]\label{lem:integrality}
		Let \(S\ge1\) be an integer, and let \(A,v\in\Q[[q]]\) satisfy \eqref{eq:integral-system}. Write \(v(q)=\sum_{n\ge0}\frac{v_n}{n!}q^n.\) Suppose
		\(c_0\in\mathbb{Z}[1/S]^\times\), \(c_1,c_2,v_1\in\mathbb{Z}[1/S]\), \(v_0=0\), 
		\(G\in\mathbb{Z}[1/S][v]\), \(F=\frac{P}{Q}\), where \(P,Q\in\mathbb{Z}[1/S][v]\) and
		\(Q(0)\in\mathbb{Z}[1/S]^\times\). Then
		\(c_n,v_n\in\mathbb{Z}[1/S]\) for every \(n\ge0\).
	\end{lemma}
	\begin{proof}
		Adjoin \(a=A^{-1}\) and \(b=Q(v)^{-1}\), which satisfy
		\[
		a'=-A'a^2,\qquad b'=-Q'(v)v'b^2.
		\]
		Here \(a(0)=c_0^{-1}\) and \(b(0)=Q(0)^{-1}\) lie in \(\mathbb{Z}[1/S]\) by the unit hypotheses.
		Together with \(A,A',A'',v,v'\), these variables satisfy a first-order
		polynomial differential system over \(\mathbb{Z}[1/S]\), with initial values in \(\mathbb{Z}[1/S]\).
		Repeated differentiation and the product rule give \(c_n,v_n\in\mathbb{Z}[1/S]\) for every \(n\ge0\).
	\end{proof}
	
	\begin{table}[!htbp]
		\centering\small
		\caption{Coordinate changes $v=\phi(u)$.}
		\label{tab:integrality-data}
		\begin{tabular}{@{}lc@{\qquad\qquad}lc@{}}
			\toprule
			\(f\)&\(v=\phi(u)\)&\(f\)&\(v=\phi(u)\)\\
			\midrule
			\(14a\)&\(\dfrac{u-16}{8}\)&\(26a\)&\(\dfrac{u-1}{4}\)\\\addlinespace[3pt]
			\(15a\)&\(\dfrac{u-9}{3}\)&\(26b\)&\(\dfrac{u-1}{4}\)\\\addlinespace[3pt]
			\(21a\)&\(\dfrac{441-25u}{16u}\)&\(38a\)&\(\dfrac{u-1}{4}\)\\\addlinespace[3pt]
			\(33a\)&\(\dfrac{9-25u}{16u}\)&\(38b\)&\(\dfrac{u-1}{4}\)\\\addlinespace[3pt]
			\(34a\)&\(\dfrac{5(4u+15)}{11(15u-16)}\)&\(58a\)&\(\dfrac{u+11}{176}\)\\\addlinespace[3pt]
			\(46a\)&\(81u-3136\)&\(58b\)&\(\dfrac{11u+1}{16}\)\\
			\bottomrule
		\end{tabular}
	\end{table}
	
	\begin{proof}[Proof of \cref{prop:Z[1/S]}]
		Substitution of the data in \cref{tab:main-data} and the coordinate changes in
		\cref{tab:integrality-data} into \eqref{eq:integral-FG} gives
		\(P,Q\in\mathbb{Z}[1/S_f][v]\) and \(G\in\mathbb{Z}[1/S_f][v]\).
		For each \(f\), the Sage certificate records \(G(v)\), \(P(v)\), \(Q(v)\), and \(Q(0)\).
		Moreover, \(Q(0) \in \mathbb{Z}[1/S_f]^\times\),
		and \(c_2\) and \(v_1=\phi'(u_0)u_1\) lie in \(\mathbb{Z}[1/S_f]\).
		Thus \cref{lem:integrality} gives \(c_n\in\mathbb{Z}[1/S_f]\) for all \(n\).
		Finally, the minimality of $S_f$ follows from the
		denominator of \(c_2\) or \(c_3\) in \cref{tab:cm-expansions}.
	\end{proof}
	
	\section{Arithmetic applications}
	\label{sec:application}
	We use the Taylor expansion of $f$ to compute the Shimura uniformization
	\[\Phi:\operatorname{Jac}(X_N) \to E,\]
	to high accuracy for the elliptic curves $E$ associated with the forms $f$ in \cref{tab:main-data}.
	
	Following \cite[Sections~1.2 and~4]{Nelson}, we consider the map
	\[
	\operatorname{Pic}^0 X_B(\mathbb{C}) \cong \operatorname{Jac} X_B(\mathbb{C}) \to E(\mathbb{C}) \cong \mathbb{C}/\Lambda
	\]
	Up to isogeny, \(\Phi\) agrees with the linear extension of
	\[
	[\tau_1] - [\tau_0] \mapsto \int_{\tau_0}^{\tau_1} f_N(z)\, dz=\int_{w(\tau_0)}^{w(\tau_1)}(\tau-\bar\tau) \frac{f_N(w)}{(1-w)^2} \, dw.
	\]
	
	Since the recurrence determines $f_N(w)$, we can compute \(\Phi\) to high accuracy. We compute the images of the CM divisors $[w\tau_D]-[\tau_D]$, where $w$ is an Atkin--Lehner operator and $\tau_D$ is a CM point of discriminant $D$, for all $w \in W_N$ and all fundamental discriminants $D$ with $|D|<200$. We present a few examples \cref{tab:further-cm-images}. The complete data are attached in the certificate. 
	\begin{table}[!htbp]
		\centering\footnotesize
		\renewcommand{\arraystretch}{1.2}
		\setlength{\tabcolsep}{2pt}
		\caption{Numerically recognized images of CM divisors on \(E\) (Cremona labels), with \(P_{D,w}=\Phi([w\tau_D]-[\tau_D])\).}
		\label{tab:further-cm-images}
		\begin{tabular}{|c|r|c|c|l|}
			\hline
			\(E\)&\(D\)&\(h(D)\)&\(w\)&Minimal polynomial of \(x(P_{D,w})\)\\
			\hline
			\(14a1\)&\(-163\)&\(1\)&\(w_{2}\)&\(150338180418782421668964075x+144420233907863660158612249\)\\
			\hline
			\(15a1\)&\(-127\)&\(5\)&\(w_{3}\)&\(\begin{aligned}
				&5166015625x^{5}+11149788901250x^{4}-84764507545775x^{3}\\
				&{}+305940630600433x^{2}-590593001260921x+563606776197313
			\end{aligned}\)\\
			\hline
			\(58a1\)&\(-131\)&\(5\)&\(w_{29}\)&\(83521x^{5}-156927x^{4}+139246x^{3}-185426x^{2}+211985x-88687\)\\
			\hline
		\end{tabular}
	\end{table}

	These minimal polynomials were recognized using the LLL algorithm~\cite[Section~1]{LLL1982}. Confidently recognizing all points within range $|D|<200$ requires approximately $4400$ coefficients. Computing this many without the recurrences is impractical.
	
	\appendix
	\section{Assumed initial values and complete CM expansions}
	\label{sec:cm-expansion}
	This appendix records which initial coefficients $c_n$ are assumed for \cref{thm:main} and gives the complete CM expansions of the eigenforms. \Cref{tab:input-data} lists the assumed coefficients for each form; their values are displayed in the expansions of \cref{tab:cm-expansions}. These values are obtained by numerical methods \cite{Nelson} within range of $10^{-40}$ and recognized as rationals by LLL~\cite[Section~1]{LLL1982}. It would be interesting to prove the values of $c_n$ using geometric methods, which makes the recurrence and Taylor expansion unconditional. 
	
	\begingroup
	\small
	\renewcommand{\arraystretch}{1.45}
	\begin{longtable}{@{}ll@{}}
		\caption{Assumed initial coefficients $c_n$ needed.}
		\label{tab:input-data}\\
		\toprule
		Forms&Assumed coefficients\\
		\midrule
		\endfirsthead
		\bottomrule
		\endlastfoot
		\(14a,\ 15a,\ 21a\)&\(c_2,\ c_3\)\\
		\(33a,\ 34a,\ 46a,\ 26a,\ 26b,\ 38a,\ 38b,\ 58b\)&\(c_2,\ c_3,\ c_4,\ c_5\)\\
		\(58a\)&\(c_2,\ c_3,\ c_4,\ c_5,\ c_6\)\\
	\end{longtable}
	\endgroup
	
	Recall the expansion \eqref{eq:f(thetaw)}:
	\[
	f_N(z)=(1-w_N)^2\sum_{n\ge0}
	c_{N,n}\frac{1}{n!}(\Theta_Nw_N)^n.
	\]
	In \cref{tab:cm-expansions}, we record all the Taylor expansions of eigenforms in \cref{tab:main-data}. In each row, \(\Theta_N=\frac{b_1}{b_0}\). The displayed algebraic expressions for
	\(\Theta_N/(\pi\Omega_D^2)\) are numerically recognized within range of $10^{-50}$. Apart from \(c_0=c_1=1\), the displayed coefficients are the assumed initial values of \cref{tab:input-data}, except the last one in each expansion, marked with \(\dagger\), which is computed from the recurrence \eqref{eq:coupled-recurrence}.
	
	\begingroup
	\small
	\renewcommand{\arraystretch}{1.35}
	\setlength{\tabcolsep}{4pt}
	\begin{longtable}{@{}rcl@{}}
		\caption{Explicit CM expansions for all twelve eigenforms. The last coefficient marked with \(\dagger\) is computed from the recurrence \eqref{eq:coupled-recurrence}; the other coefficients \(c_n\) with \(n\ge2\) are the assumed initial values.}
		\label{tab:cm-expansions}\\
		\toprule
		Form&\(\tau_N\)&\(\Theta_N\)\\
		\midrule
		\endfirsthead
		\multicolumn{3}{c}{\tablename~\thetable\ (continued)}\\
		\toprule
		Form&\(\tau_N\)&\(\Theta_N\)\\
		\midrule
		\endhead
		\midrule
		\multicolumn{3}{r}{\footnotesize Continued on the next page}\\
		\endfoot
		\bottomrule
		\endlastfoot
		14a&\(\dfrac{-\sqrt{14}+\sqrt{-22}}6\)
		&\(\pi\Omega_{-11}^2\left(1+\dfrac17\sqrt{-77}\right)\)\\*[2pt]
		\multicolumn{3}{@{}l@{}}{\(\displaystyle\begin{aligned}
				f_{14}(z)=(1-w_{14})^2\Bigl(&1+\Theta_{14}w_{14}-4\frac{1}{2!}(\Theta_{14}w_{14})^2 {}+\frac{214}{9}\frac{1}{3!}(\Theta_{14}w_{14})^3\\
				&{}+24^{\dagger}\frac{1}{4!}(\Theta_{14}w_{14})^4+O(w_{14}^5)\Bigr).
			\end{aligned}\)}\\[4pt]
		\midrule
		15a&\(\dfrac{\sqrt{-210}-\sqrt{-14}}{14}\)
		&\(-\dfrac8{15}\pi\Omega_{-7}^2\sqrt{15}\)\\*[2pt]
		\multicolumn{3}{@{}l@{}}{\(\displaystyle\begin{aligned}
				f_{15}(z)=(1-w_{15})^2\Bigl(&1+\Theta_{15}w_{15}-\frac{1}{4}\frac{1}{2!}(\Theta_{15}w_{15})^2 {}+\frac{27}{16}\frac{1}{3!}(\Theta_{15}w_{15})^3\\
				&{}-\frac{381}{64}^{\dagger}\frac{1}{4!}(\Theta_{15}w_{15})^4+O(w_{15}^5)\Bigr).
			\end{aligned}\)}\\[4pt]
		\midrule
		21a&\(\dfrac{-21-5\sqrt{21}+5\sqrt{-43}+\sqrt{-903}}{16}\)
		&\(\pi\Omega_{-43}^2\left(-\dfrac{63}{4}-\dfrac34\sqrt{-903}\right)\)\\*[2pt]
		\multicolumn{3}{@{}l@{}}{\(\displaystyle\begin{aligned}
				f_{21}(z)=(1-w_{21})^2\Bigl(&1+\Theta_{21}w_{21}-\frac{905}{3969}\frac{1}{2!}(\Theta_{21}w_{21})^2 {}-\frac{1318}{27783}\frac{1}{3!}(\Theta_{21}w_{21})^3\\
				&{}-\frac{4357004}{5250987}^{\dagger}\frac{1}{4!}(\Theta_{21}w_{21})^4+O(w_{21}^5)\Bigr).
			\end{aligned}\)}\\[4pt]
		\midrule
		33a&\(\dfrac{-\sqrt{33}+\sqrt{-67}}{10}\)
		&\(\pi\Omega_{-67}^2\left(-\dfrac{23}{5}-\dfrac{23}{165}\sqrt{-2211}\right)\)\\*[2pt]
		\multicolumn{3}{@{}l@{}}{\(\displaystyle\begin{aligned}
				f_{33}(z)=(1-w_{33})^2\Bigl(&1+\Theta_{33}w_{33}-\frac{72337}{529}\frac{1}{2!}(\Theta_{33}w_{33})^2 {}+\frac{9202302}{12167}\frac{1}{3!}(\Theta_{33}w_{33})^3\\
				&{}+\frac{3009372468}{279841}\frac{1}{4!}(\Theta_{33}w_{33})^4 {}-\frac{1149714703920}{6436343}\frac{1}{5!}(\Theta_{33}w_{33})^5\\
				&{}+\frac{176790432531276}{148035889}^{\dagger}\frac{1}{6!}(\Theta_{33}w_{33})^6+O(w_{33}^7)\Bigr).
			\end{aligned}\)}\\[4pt]
		\midrule
		34a&\(\begin{gathered}
			\bigl(-34\sqrt3+7\sqrt{102}\\[-2pt]
			{}+7\sqrt{-33}-\sqrt{-1122}\bigr)/45
		\end{gathered}\)
		&\(\pi\Omega_{-11}^2\left(\begin{gathered}
			1-\dfrac{11}{102}\sqrt{34}\\[-2pt]
			{}-\dfrac13\sqrt{-11}-\dfrac1{34}\sqrt{-374}
		\end{gathered}\right)\)\\*[2pt]
		\multicolumn{3}{@{}l@{}}{\(\displaystyle\begin{aligned}
				f_{34}(z)=(1-w_{34})^2\Bigl(&1+\Theta_{34}w_{34}-\frac{12}{25}\frac{1}{2!}(\Theta_{34}w_{34})^2 {}+\frac{2762}{125}\frac{1}{3!}(\Theta_{34}w_{34})^3\\
				&{}-\frac{17736}{625}\frac{1}{4!}(\Theta_{34}w_{34})^4 {}-\frac{67648}{125}\frac{1}{5!}(\Theta_{34}w_{34})^5\\
				&{}+\frac{6553872}{15625}^{\dagger}\frac{1}{6!}(\Theta_{34}w_{34})^6+O(w_{34}^7)\Bigr).
			\end{aligned}\)}\\[4pt]
		\midrule
		46a&\(\dfrac{35\sqrt{-489}+4\sqrt{-22494}}{489}\)
		&\(\dfrac{7855}{69}\pi\Omega_{-163}^2\sqrt{46}\)\\*
		\addlinespace[2.5pt]
		\multicolumn{3}{@{}l@{}}{\(\displaystyle\begin{aligned}
				f_{46}(z)&=(1-w_{46})^2\Bigl(1+\Theta_{46}w_{46}+\frac{27824049}{123402050}\frac{1}{2!}(\Theta_{46}w_{46})^2 {}-\frac{114816966372}{484661551375}\frac{1}{3!}(\Theta_{46}w_{46})^3\\
				&\qquad{}-\frac{784666539922896}{3807016486050625}\frac{1}{4!}(\Theta_{46}w_{46})^4 {}-\frac{3694325348311426944}{5980822899585531875}\frac{1}{5!}(\Theta_{46}w_{46})^5\\
				&\qquad{}+\frac{939622156860406531404672}{234896819381221764390625}^{\dagger}\frac{1}{6!}(\Theta_{46}w_{46})^6+O(w_{46}^7)\Bigr).
			\end{aligned}\)}\\[4pt]
		\midrule
		26a&\(\begin{gathered}\bigl(52\sqrt5+14\sqrt{130}\\[-2pt]{}+7\sqrt{-55}+\sqrt{-1430}\bigr)/115\end{gathered}\)
		&\(-\dfrac{21}{2990}\pi\Omega_{-11}^2\left(\begin{gathered}988-33\sqrt{26}\\[-2pt]{}-156\sqrt{-11}-19\sqrt{-286}\end{gathered}\right)\)\\*[2pt]
		\multicolumn{3}{@{}l@{}}{\(\displaystyle\begin{aligned}
				f_{26,a}(z)=(1-w_{26})^2\Bigl(&1+\Theta_{26}w_{26}-\frac{16}{147}\frac{1}{2!}(\Theta_{26}w_{26})^2 {}+\frac{274}{3087}\frac{1}{3!}(\Theta_{26}w_{26})^3\\
				&{}-\frac{2488}{9261}\frac{1}{4!}(\Theta_{26}w_{26})^4 {}+\frac{352768}{1361367}\frac{1}{5!}(\Theta_{26}w_{26})^5\\
				&{}+\frac{5961776}{28588707}^{\dagger}\frac{1}{6!}(\Theta_{26}w_{26})^6+O(w_{26}^7)\Bigr).
			\end{aligned}\)}\\[4pt]
		\midrule
		26b&\(\begin{gathered}\bigl(52\sqrt5+14\sqrt{130}\\[-2pt]{}+7\sqrt{-55}+\sqrt{-1430}\bigr)/115\end{gathered}\)
		&\(\dfrac{1}{2990}\pi\Omega_{-11}^2\left(\begin{gathered}988-33\sqrt{26}\\[-2pt]{}-156\sqrt{-11}-19\sqrt{-286}\end{gathered}\right)\)\\*[2pt]
		\multicolumn{3}{@{}l@{}}{\(\displaystyle\begin{aligned}
				f_{26,b}(z)=(1-w_{26})^2\Bigl(&1+\Theta_{26}w_{26}+40\frac{1}{2!}(\Theta_{26}w_{26})^2 {}-734\frac{1}{3!}(\Theta_{26}w_{26})^3\\
				&{}-11240\frac{1}{4!}(\Theta_{26}w_{26})^4 {}+1308544\frac{1}{5!}(\Theta_{26}w_{26})^5\\
				&{}+5134480^{\dagger}\frac{1}{6!}(\Theta_{26}w_{26})^6+O(w_{26}^7)\Bigr).
			\end{aligned}\)}\\[4pt]
		\midrule
		38a&\(\dfrac{-2\sqrt{38}+\sqrt{-209}}{19}\)
		&\(-\dfrac9{19}\pi\Omega_{-11}^2\left(4+\sqrt{-22}\right)\)\\*[2pt]
		\multicolumn{3}{@{}l@{}}{\(\displaystyle\begin{aligned}
				f_{38,a}(z)=(1-w_{38})^2\Bigl(&1+\Theta_{38}w_{38}-\frac{64}{27}\frac{1}{2!}(\Theta_{38}w_{38})^2 {}-\frac{62}{27}\frac{1}{3!}(\Theta_{38}w_{38})^3\\
				&{}-\frac{19672}{2187}\frac{1}{4!}(\Theta_{38}w_{38})^4 {}+\frac{41984}{6561}\frac{1}{5!}(\Theta_{38}w_{38})^5\\
				&{}-\frac{38071696}{177147}^{\dagger}\frac{1}{6!}(\Theta_{38}w_{38})^6+O(w_{38}^7)\Bigr).
			\end{aligned}\)}\\[4pt]
		\midrule
		38b&\(\dfrac{-2\sqrt{38}+\sqrt{-209}}{19}\)
		&\(\dfrac{13}{19}\pi\Omega_{-11}^2\left(4+\sqrt{-22}\right)\)\\*[2pt]
		\multicolumn{3}{@{}l@{}}{\(\displaystyle\begin{aligned}
				f_{38,b}(z)=(1-w_{38})^2\Bigl(&1+\Theta_{38}w_{38}-\frac{16}{169}\frac{1}{2!}(\Theta_{38}w_{38})^2 {}-\frac{878}{2197}\frac{1}{3!}(\Theta_{38}w_{38})^3\\
				&{}+\frac{7688}{28561}\frac{1}{4!}(\Theta_{38}w_{38})^4 {}-\frac{2699648}{371293}\frac{1}{5!}(\Theta_{38}w_{38})^5\\
				&{}+\frac{25568336}{4826809}^{\dagger}\frac{1}{6!}(\Theta_{38}w_{38})^6+O(w_{38}^7)\Bigr).
			\end{aligned}\)}\\[4pt]
		\midrule
		58a&\(\dfrac{-3\sqrt{58}+\sqrt{-319}}{29}\)
		&\(-\dfrac5{29}\pi\Omega_{-11}^2\left(6+\sqrt{-22}\right)\)\\*[2pt]
		\multicolumn{3}{@{}l@{}}{\(\displaystyle\begin{aligned}
				f_{58,a}(z)=(1-w_{58})^2\Bigl(&1+\Theta_{58}w_{58}-\frac{4}{5}\frac{1}{2!}(\Theta_{58}w_{58})^2 {}+\frac{314}{125}\frac{1}{3!}(\Theta_{58}w_{58})^3\\
				&{}+\frac{31192}{625}\frac{1}{4!}(\Theta_{58}w_{58})^4 {}-\frac{1977728}{625}\frac{1}{5!}(\Theta_{58}w_{58})^5\\
				&{}+\frac{4428752}{3125}\frac{1}{6!}(\Theta_{58}w_{58})^6 {}-\frac{4829523664}{78125}^{\dagger}\frac{1}{7!}(\Theta_{58}w_{58})^7+O(w_{58}^8)\Bigr).
			\end{aligned}\)}\\[4pt]
		\midrule
		58b&\(\dfrac{-3\sqrt{58}+\sqrt{-319}}{29}\)
		&\(\dfrac3{29}\pi\Omega_{-11}^2\left(6+\sqrt{-22}\right)\)\\*[2pt]
		\multicolumn{3}{@{}l@{}}{\(\displaystyle\begin{aligned}
				f_{58,b}(z)=(1-w_{58})^2\Bigl(&1+\Theta_{58}w_{58}+\frac{140}{9}\frac{1}{2!}(\Theta_{58}w_{58})^2 {}+\frac{262}{3}\frac{1}{3!}(\Theta_{58}w_{58})^3\\
				&{}-\frac{39592}{81}\frac{1}{4!}(\Theta_{58}w_{58})^4 {}+\frac{420224}{81}\frac{1}{5!}(\Theta_{58}w_{58})^5\\
				&{}-\frac{303950320}{729}^{\dagger}\frac{1}{6!}(\Theta_{58}w_{58})^6+O(w_{58}^7)\Bigr).
			\end{aligned}\)}\\[4pt]
	\end{longtable}
	\endgroup
	
	\section{Certification of CM coordinates}
	\label{app:u-values}
	We verify the values of the Hauptmodul $u$ at the CM points on $X_N$ used in \cref{tab:main-data} and \cref{tab:families}.
	
	First, for \(N=15,21,33\), the models in \cref{tab:models} are preserved by
	\((x,y)\mapsto(m_N/x,m_Ny/x^2)\), with
	\(()m_{15},m_{21},m_{33})=(3,7,9)\)
	from Gonz\'alez--Rotger \cite[Theorem~3.4, Table~1]{GR}. We fix this ambiguity by setting
	\[
	u(\mathrm{CM}_{15}(-3))=\infty,\qquad
	u(\mathrm{CM}_{21}(-7))=0,\qquad
	u(\mathrm{CM}_{33}(-3))=0.
	\]
	Here \(\mathrm{CM}_N(D)\) denotes the set of CM points of discriminant \(D\) on $X_N$. 
	
	Most of the $u$-values have already been verified in \cite{Elkies,GuoYang2017,LinYang2020}.
	
	\begingroup
	\small
	\renewcommand{\arraystretch}{1.1}
	\setlength{\tabcolsep}{7pt}
	\begin{longtable}{@{}rlll@{}}
		\caption{Certified \(u\)-values. For genus two, \(u=x^2\).}\\
		\toprule
		\(N\)&\(D\)&\(u\)&Certification\\
		\midrule
		\endfirsthead
		\multicolumn{4}{c}{\tablename~\thetable\ (continued)}\\
		\toprule
		\(N\)&\(D\)&\(u\)&Certification\\
		\midrule
		\endhead
		\bottomrule
		\endlastfoot
		14&\(-11,\;-4\)&\(16,\;\infty\)&
		\cite[Section~5.1]{Elkies}\\
		15&\(-7,\;-3\)&\(9,\;\infty\)&
		\cite[Table~6]{Elkies}\\
		21&\(-43\)&\(441/25\)&
		\eqref{eq:cm-extra-certification} below\\
		21&\(-4\)&\(49\)&
		\cite[Table~6]{LinYang2020}\\
		33&\(-67\)&\(9/25\)&
		\eqref{eq:cm-extra-certification} below\\
		33&\(-3,\;-4\)&\(0,\;1\)&
		\cite[Table~6]{LinYang2020}\\
		34&\(-11,\;-3\)&\(-15/4,\;16/15\)&
		\cite[Table~6]{LinYang2020}\\
		46&\(-163\)&\(3136/81\)&
		\eqref{eq:cm-extra-certification} below\\
		46&\(-4,\;-3\)&\(\infty,\;64\)&
		\cite[Table~6]{LinYang2020}\\
		\midrule
		26&\(-11,\;-52,\;-8\)&\(1,\;0,\;\infty\)&
		\cite[Example~33]{GuoYang2017}\\
		38&\(-11,\;-19,\;-4\)&\(1,\;0,\;\infty\)&
		\cite[Table~11]{GuoYang2017}\\
		58&\(-11,\;-232,\;-8,\;-3\)&\(-11,\;0,\;\infty,\;-27\)&
		\cite[Table~16]{GuoYang2017}\\
	\end{longtable}
	\endgroup
	For Elkies's coordinates \(T_N\) \cite[equations~(68), (76)]{Elkies} at $N=14$ and $15$, the changes of variable are
	\[u_{14}=-16T_{14},\qquad u_{15}=-T_{15}/3\]
	
	For Lin--Yang's \(j_N\) \cite[Remark~55, Table~6 and Table~7]{LinYang2020}, the coordinate changes are
	\[
	x=7/j_{21},\qquad x=1/j_{33},\qquad
	x=\frac{5+3j_{34}}{3-5j_{34}},\qquad x=8j_{46},
	\]

	The following three values are not directly recorded in the literature.
	\[
	u(\mathrm{CM}_{21}(-43))=\frac{441}{25},\qquad
	u(\mathrm{CM}_{33}(-67))=\frac9{25},\qquad
	u(\mathrm{CM}_{46}(-163))=\frac{3136}{81}.
	\]
	The discriminants $-43$, $-67$, and $-163$ are also admissible on $X_6$. We therefore use Lin--Yang's parametrization $S_N$ of the associated principally polarized abelian surfaces in $\mathcal{A}_2$ by $X_N$ \cite[equation~(15), Example~29 and Table~4]{LinYang2020}. For \(N=21,33,46\), Lin--Yang give explicit maps
	\[
	\begin{aligned}
		S_N:X_N/\langle w_N\rangle&\longrightarrow\mathcal A_2,\\
		j&\longmapsto[f_2(j):f_3(j):f_5(j):f_6(j)],
	\end{aligned}
	\]
	For discriminant \(6\), Lin--Yang \cite[Example~29]{LinYang2020} give
	\[
	\begin{aligned}
		S_6:X_6/W_6&\longrightarrow\mathcal A_2,\\
		J&\longmapsto[J:J^2:J^3:J^3(J+1)],
	\end{aligned}
	\]
	with \(J\) taking the values \(0,-16/27,\infty\) at CM discriminants
	\(-4,-24,-3\), respectively

	Direct substitution in Lin--Yang's parametrizations gives
	\begin{equation}\label{eq:cm-extra-moduli}
		S_N(j_D)=S_6(J_D),
	\end{equation}
	for the following candidate pairs \((j_D,J_D)\):
	\[
	\begin{array}{c|r|c|r}
		N&D&j_D&J_D\\ \hline
		21&-43&5/3&21^4/10^6\\
		33&-67&5/3&231^4/20^6\\
		46&-163&7/9&27531^4/1870^6
	\end{array}.
	\]
	
	Errthum's coordinate \(t_6\) is related to \(J_6\) by
	\(J_6=-\frac{16}{27}t_6\). Hence \cite[Table~3]{Errthum2011}
	verifies that each \(J_D\) equals \(J_6(P_D')\) for some
	\(P_D'\in\mathrm{CM}_6(D)\).
	
	To certify that \(j_D\) is a CM value, choose \(P_D\in X_N\) with
	\(j_N(P_D)=j_D\). Equation~\eqref{eq:cm-extra-moduli} identifies the
	principally polarized abelian surfaces associated with \(P_D\) and \(P_D'\).
	Denote this surface by \(A\), and let
	\(\iota_M:\mathcal O_M\hookrightarrow\operatorname{End}(A)\)
	be the actions of the corresponding maximal quaternion orders.
	Write \(\operatorname{End}(A,\iota_M)\) for the ring of endomorphisms
	commuting with \(\iota_M\).
	Since $P_D'\in\mathrm{CM}_6(D)$, the moduli interpretation gives
	$\operatorname{End}(A,\iota_6)\simeq\mathcal O_D$
	\cite[\S2.2]{Molina2012}.
	
	The complex span of either quaternionic action on
	$\operatorname{Lie}(A)$ is $M_2(\mathbb{C})$, whose commutant
	consists of scalar maps
	\cite[proof of Main Theorem~43.6.14]{Voight2021}.
	Since an endomorphism of a complex abelian variety is determined
	by its differential \cite[proof of Proposition~43.6.28]{Voight2021}, we obtain
	\[
	\operatorname{End}(A,\iota_N)
	=
	\{f\in\operatorname{End}(A):df\text{ is scalar}\}
	=
	\operatorname{End}(A,\iota_6)
	\simeq\mathcal O_D.
	\]
	Hence $P_D\in\mathrm{CM}_N(D)$.
	We have thus identified CM points
	\(
	P_{21}\in\mathrm{CM}_{21}(-43),
	P_{33}\in\mathrm{CM}_{33}(-67),\) and
	\(P_{46}\in\mathrm{CM}_{46}(-163)
	\)
	satisfying $j_{21}(P_{21})=\frac53$, $	j_{33}(P_{33})=\frac53$, and $j_{46}(P_{46})=\frac79$.
	
	The coordinate changes above therefore give
	\begin{equation}\label{eq:cm-extra-certification}
		\begin{aligned}
			u(\mathrm{CM}_{21}(-43))
			&=\frac{49}{(5/3)^2}=\frac{441}{25},\\
			u(\mathrm{CM}_{33}(-67))
			&=\frac{1}{(5/3)^2}=\frac9{25},\\
			u(\mathrm{CM}_{46}(-163))
			&=8^2(7/9)^2=\frac{3136}{81}.
		\end{aligned}
	\end{equation}
	
	\section*{Computational certificates}
	All computational certificates for this paper are available in
	\href{https://github.com/jxgeng-2/shimura_recurrence_certificates/releases/tag/v1.1.0}
	{version~1.1.0 of the accompanying GitHub repository}.
	
	\section*{AI disclosure}
	AI tools were extensively used throughout the research for this paper, including for
	exploring the recursion, suggesting ideas, writing code, examining proofs, and improving formatting.
	The paper is human-written with assistance of AI and the author takes full responsibility for the correctness of the results.

	\begin{center}
		\begin{tabular}{c}
			\textbf{Jingxuan Geng} \\
			Department of Mathematics, The Ohio State University \\
			\texttt{geng.210@osu.edu}
		\end{tabular}
	\end{center}
\end{document}